\documentstyle[11pt]{article}

\begin{document}
\noindent
\begin{center}
  {\LARGE A New Cohomology Theory of Orbifold}\footnote{both authors partially
supported by the National Science Foundation}
  \end{center}

  \noindent
  \begin{center}

    {\large Weimin Chen\footnote{current address: SUNY-Stony Brook,
NY 11794. e-mail: wechen@math.sunysb.edu} and Yongbin Ruan}\\[5pt]
      Department of Mathematics, University of Wisconsin-Madison\\
        Madison, WI 53706\\[5pt]

              \end{center}

              \def \x{{\bf x}}
              \def \M{{\cal M}}
              \def \C{{\bf C}}
              \def \Z{{\bf Z}}
              \def \R{{\bf R}}
              \def \Q{{\bf C}}
              \def \U{{\cal U}}
              \def \E{{\cal E}}
              \def \z{{\bf z}}
              \def \m{{\bf m}}
              \def \n{{\bf n}}
              \def \g{{\bf g}}
              \def \h{{\bf h}}
              \def \V{{\cal V}}
              \def \W{{\cal W}}
              \def \T{{\cal T}}
              \def \P{{\bf P}}
              \def \F{{\bf C}}

\tableofcontents

\section{Introduction}
      An orbifold is a topological space locally modeled
      on the quotient of a smooth manifold by a
      finite group. Therefore, orbifolds belong to one of the simplest kinds
      of singular
      spaces. Orbifolds appear naturally in many branches of
      mathematics.  For example, symplectic reduction often gives rise to
      orbifolds.
      An algebraic 3-fold with terminal singularities can be
      deformed into a symplectic orbifold. Orbifold also appears
      naturally in string theory, where many known Calabi-Yau
      3-folds are the so called crepant resolutions of a Calabi-Yau
      orbifold. The physicists even attempted
      to formulate string theories on Calabi-Yau orbifolds which are expected
      to be ``equivalent'' to the string
      theories on its crepant resolutions \cite{DHVW}.
      As a consequence of this orbifold string theory consideration, one has
      the following prediction that
      {\em ``orbifold quantum cohomology'' is ``isomorphic''
      to the ordinary quantum cohomology of its crepant resolutions.}
      At this moment, even the
      physical idea around this subject is still vague and incomplete,
      particularly for the possible isomorphism. However,
      it seems that there are interesting new mathematical structures
      that are behind such orbifold string theories.

      This article is the first paper of a program
      to understand these new mathematical treasures
      behind orbifold string theory. We introduce {\it orbifold cohomology
      groups} of an almost complex orbifold, and {\it orbifold
      Dolbeault cohomology groups} of a complex orbifold. The main
      result of this paper is the construction of orbifold cup products on
      orbifold cohomology groups and orbifold Dolbeault cohomology groups,
      which make the corresponding total orbifold cohomology into a ring
      with unit. We will call the resulting rings {\it orbifold cohomology
      ring} or {\it orbifold Dolbeault cohomology ring}.
(See Theorems 4.1.5 and 4.1.7 for details). In the case when the
almost complex orbifold is closed and symplectic, the orbifold cohomology
ring corresponds to the ``classical part'' of the orbifold quantum
cohomology ring constructed in \cite{CR}. Originally, this article
is a small part
      of the much longer paper \cite{CR} regarding the theory of orbifold
      quantum cohomology.
      However, we feel that the classical part (i.e. the orbifold cohomology)
      of the orbifold quantum cohomology is interesting in its own right,
      and technically, it is also much simpler to construct.
      Therefore, we decided to put it in a separate paper.

      A brief history is in order. In the case of Gorenstein global quotients,
      orbifold Euler characteristic-Hodge numbers
      have been extensively studied in the literature
      (see \cite{RO},\cite{BD},\cite{Re} for
      a more complete reference). However, we would like to point out that
      (i) our orbifold cohomology is well-defined  for any almost complex
      orbifold
      which may or may not be Gorenstein.
      Furthermore, it has an interesting feature that an orbifold cohomology
      class
      of a non-Gorenstein orbifold could have
      a rational degree (See examples in section 5); (ii) Even in the case of
      Gorensterin orbifolds, orbifold
    cohomology ring contains much more information than just orbifold
    Betti-Hodge numbers.
      In the case of global quotients, some constructions of this paper
    are already
    known to physicists. A notable exception is the orbifold cup product.
    On the other hand, many interesting orbifolds are not global quotients
    in general.  For examples, most of Calabi-Yau hypersurfaces
    of weighted projective spaces are not global quotients. In this article,
    we systematically developed the theory (including the construction of
    orbifold cup products) for general orbifolds. Our construction of
    orbifold cup products is motivated
    by the construction of orbifold quantum cohomology.

     The second author would like to thank R. Dijkgraaf for
     bringing the orbifold string theory to his attention and E.
     Zaslow for valuable discussions.

\section{Recollections on Orbifold}

In this section, we review basic definitions in the theory of
orbifold. A systematic treatment of various aspects of
differential geometry on orbifolds is contained in our forthcoming
paper \cite{CR}. The notion of orbifold was first introduced by
Satake in \cite{S}, where a different name, {\it V-manifold}, was
used.

\vspace{2mm}

Let $U$ be a connected topological space, $V$ be a connected
n-dimensional smooth manifold with a smooth action by a finite
group $G$. Here we assume throughout that {\it the fixed-point set of
each element of the group is either the whole space  or of codimension at
least two}.
In particular, the action of $G$ does not have to be effective. This is the
case, for example, when the action is orientation-preserving. This
requirement has a consequence that the non-fixed-point set is
locally connected if it is not empty. We will call the subgroup of $G$, which
consists of elements fixing the whole space $V$, the {\it kernel} of the
action. {\it An n-dimensional uniformizing system} of $U$ is a triple
$(V,G,\pi)$, where $\pi: V \rightarrow U$ is a continuous map
inducing a homeomorphism between the quotient space $V/G$ and $U$.
Two uniformizing systems $(V_i,G_i,\pi_i)$, $i=1,2$, are {\it
isomorphic} if there is a diffeomorphism $\phi: V_1\rightarrow
V_2$ and an isomorphism $\lambda: G_1\rightarrow G_2$ such that
$\phi$ is $\lambda$-equivariant, and $\pi_2\circ\phi=\pi_1$. If
$(\phi,\lambda)$ is an automorphism of $(V,G,\pi)$, then
there is  $g\in G$ such that $\phi(x)=g\cdot x$ and
$\lambda(a)=gag^{-1}$ for any $x\in V$ and $a\in G$. Note that here $g$ is
unique iff the action of $G$ on $V$ is effective.

Let $i:U^\prime\rightarrow U$ be a connected open subset of $U$,
and $(V^\prime,G^\prime,\pi^\prime)$ be a uniformizing system of
$U^\prime$. We say that $(V^\prime,G^\prime,\pi^\prime)$ is
induced from $(V,G,\pi)$ if there is a monomorphism $\tau:
G^\prime\rightarrow G$ which is an isomorphism restricted to the kernels of
the action of $G^\prime$ and $G$ respectively, and a $\tau$-equivariant open
embedding $\psi: V^\prime\rightarrow V$ such that
$i\circ\pi^\prime=\pi\circ\psi$. The pair $(\psi,\tau):
(V^\prime,G^\prime,\pi^\prime)\rightarrow (V,G,\pi)$ is called an
{\it injection}. Two injections $(\psi_i,\tau_i):
(V^\prime_i,G^\prime_i,\pi^\prime_i)\rightarrow (V,G,\pi)$,
$i=1,2$, are {\it isomorphic} if there is an isomorphism
$(\phi,\lambda)$ between $(V^\prime_1,G^\prime_1,\pi^\prime_1)$
and $(V^\prime_2,G^\prime_2,\pi^\prime_2)$, and an automorphism
$(\bar{\phi}, \bar{\lambda})$ of $(V,G,\pi)$ such that
$(\bar{\phi},\bar{\lambda})\circ
(\psi_1,\tau_1)=(\psi_2,\tau_2)\circ (\phi,\lambda)$. One can
easily verify that for any connected open subset $U^\prime$ of
$U$, a uniformizing system $(V,G,\pi)$ of $U$ induces a unique
isomorphism class of uniformizing systems of $U^\prime$.

Let $U$ be a connected and locally connected topological space.
For any point $p\in U$, we can define the {\it germ} of
uniformizing systems at $p$ in the following sense. Let $(V_1,
G_1,\pi_1)$ and $(V_2,G_2,\pi_2)$ be uniformizing systems of
neighborhoods $U_1$ and $U_2$ of $p$. We say that $(V_1,
G_1,\pi_1)$ and $(V_2,G_2,\pi_2)$ are {\it equivalent} at $p$ if
they induce isomorphic uniformizing systems for a neighborhood
$U_3$ of $p$.

\vspace{3mm}

\noindent{\bf Definition 2.1:}{\it
\begin{enumerate}
\item Let $X$ be a Hausdorff, second countable topological space.
An {\it n-dimensional orbifold structure} on $X$ is given by the
following data: for any point $p\in X$, there is a neighborhood
$U_p$ and a n-dimensional uniformizing system $(V_p, G_p, \pi_p)$
of $U_p$ such that for any point $q\in U_p$, $(V_p, G_p, \pi_p)$
and $(V_q,G_q,\pi_q)$ are equivalent at $q$ (i.e., defining the
same germ at $q$).  With a given germ of orbifold structures, $X$
is called an {\it orbifold}.  An open subset $U$ of $X$ is called
a {\it uniformized open subset} if it is uniformized by a
$(V,G,\pi)$ such that for each $p\in U$, $(V,G,\pi)$ defines the
same germ with $(V_p,G_p,\pi_p)$ at $p$. We may assume that each $V_p$ is
a n-ball centered at origin $o$ and $\pi_p^{-1}(p)=o$. In particular,
the origin $o$ is fixed by $G_p$. If $G_p$ acts effectively
for every $p$, we call $X$ {\em a reduced orbifold.}
\item The notion of {\it orbifold with boundary}, in which
we allow the uniformizing systems to be smooth manifolds with
boundary, with a finite group action preserving the boundary, can
be similarly defined. If $X$ is an orbifold with boundary, then it
is easily seen that the boundary $\partial X$ inherits an orbifold
structure from $X$ and becomes an orbifold.
\end{enumerate}}
\hfill $\Box$

In a certain sense, Satake's definition of orbifold is less
intrinsic than ours, although they are equivalent. In \cite{S}, an
orbifold structure on $X$ is given by an open cover
$\U$ of $X$ satisfying the following conditions:
\begin{itemize}
\item [{1.}] Each element $U$ in $\U$ is uniformized, say by $(V,G,\pi)$.
\item [{2.}] If $U^\prime\subset U$, then there is a collection of
injections $(V^\prime,G^\prime,\pi^\prime)\rightarrow (V,G,\pi)$.
\item [{3.}] For any point $p\in U_1\cap U_2$, $U_1,U_2\in\U$,
there is a $U_3\in\U$ such that $p\in U_3\subset U_1\cap U_2$.
\end{itemize}
It clearly defines an orbifold structure on $X$ in the sense of Definition 2.1.
We will call such a cover of an orbifold $X$ a {\it compatible
cover} if it gives rise to the same germ of orbifold structures on $X$.
We remark that the orbifolds considered by Satake in \cite{S} are all reduced.

\vspace{2mm}

Now  we consider a class of continuous maps between two orbifolds
which respect the orbifold structures in a certain sense. Let $U$
be uniformized by $(V,G,\pi)$ and $U^\prime$ by
$(V^\prime,G^\prime,\pi^\prime)$, and $f: U\rightarrow U^\prime$
be a continuous map. A {\it $C^l$ lifting, $0\leq l\leq\infty$},
of $f$ is a $C^l$ map $\tilde{f}:V\rightarrow V^\prime$ such that
$\pi^\prime\circ\tilde{f}=f\circ\pi$, and for any $g\in G$, there
is $g^\prime \in G^\prime$ so that
$g^\prime\cdot\tilde{f}(x)=\tilde{f}(g\cdot x)$ for any $x\in V$.
Two liftings $\tilde{f}_i: (V_i,G_i,\pi_i)\rightarrow
(V^\prime_i,G^\prime_i,\pi^\prime_i)$, $i=1,2$, are {\it
isomorphic} if there exist isomorphisms
$(\phi,\tau):(V_1,G_1,\pi_1)\rightarrow (V_2,G_2,\pi_2)$ and
$(\phi^\prime,\tau^\prime):(V_1^\prime,G_1^\prime,\pi_1^\prime)\rightarrow
(V_2^\prime,G_2^\prime,\pi_2^\prime)$ such that $\phi^\prime\circ
\tilde{f}_1=\tilde{f}_2\circ \phi$.

Let $p\in U$ be any point. Then for any uniformized neighborhood
$U_p$ of $p$ and uniformized neighborhood $U_{f(p)}$ of $f(p)$
such that $f(U_p)\subset U_{f(p)}$, a lifting $\tilde{f}$ of $f$
will induce a lifting $\tilde{f}_p$ for $f|_{U_{p}}:U_p\rightarrow
U_{f(p)}$ as follows: For any injection
$(\phi,\tau):(V_p,G_p,\pi_p)\rightarrow (V,G,\pi)$, consider the
map $\tilde{f}\circ\phi: V_p\rightarrow V^\prime$. Observe that
the inclusion $\pi^\prime\circ\tilde{f}\circ\phi(V_p)\subset
U_{f(p)}$ implies that $\tilde{f}\circ\phi(V_p)$ lies in
$(\pi^\prime)^{-1}(U_{f(p)})$. Therefore there is an injection
$(\phi^\prime,\tau^\prime):(V_{f(p)}, G_{f(p)},\pi_{f(p)})
\rightarrow (V^\prime, G^\prime, \pi^\prime)$ such that
$\tilde{f}\circ\phi(V_p)\subset \phi^\prime(V_{f(p)})$. We define
$\tilde{f}_p=(\phi^\prime)^{-1}\circ\tilde{f}\circ\phi$. In this
way we obtain a lifting $\tilde{f}_p: (V_p,G_p,\pi_p)\rightarrow
(V_{f(p)}, G_{f(p)},\pi_{f(p)})$ for $f|_{U_{p}}:U_p\rightarrow
U_{f(p)}$. We can verify that different choices give isomorphic
liftings. We define the {\it germ} of liftings as follows: two
liftings are {\it equivalent at $p$} if they induce isomorphic
liftings on a smaller neighborhood of $p$.

Let $f:X\rightarrow X^\prime$ be a continuous map between
orbifolds $X$ and $X^\prime$. A {\it lifting} of $f$ consists of
following data: for any point $p \in X$, there exist charts
$(V_p,G_p,\pi_p)$ at $p$ and $(V_{f(p)},G_{f(p)}, \pi_{f(p)})$ at
$f(p)$ and a lifting $\tilde{f}_p$ of $f|_{\pi_p(V_p)}:
\pi_p(V_p)\rightarrow \pi_{f(p)}(V_{f(p)})$ such that for any
$q\in \pi_p(V_p)$, $\tilde{f}_p$ and $\tilde{f}_q$ induce the same
germ of liftings of $f$ at $q$. We can define the {\it germ} of
liftings in the sense that two liftings of $f$,
$\{\tilde{f}_{p,i}: (V_{p,i},G_{p,i},\pi_{p,i}) \rightarrow
(V_{f(p),i},G_{f(p),i},\pi_{f(p),i}): p\in X\}$, $i=1,2$, are {\it
equivalent} if for each $p\in X$, $\tilde{f}_{p,i}, i=1,2$, induce
the same germ of liftings of $f$ at $p$.

\vspace{2mm}

\noindent{\bf Definition 2.2: }{\it  A {\it $C^l$ map} ($0\leq
l\leq\infty$) between orbifolds $X$ and $X^\prime$ is a germ of
$C^l$ liftings of a continuous map between $X$ and $X^\prime$.
\hfill $\Box$}

\vspace{2mm}

We denote by $\tilde{f}$ a $C^l$ map which is a germ of liftings
of a continuous map $f$. Our definition of $C^l$ maps corresponds
to the notion of {\it $V$-maps} in [S].

\vspace{2mm}

Next we describe the notion of orbifold vector bundle, which
corresponds to the notion of smooth vector bundle over manifolds.
When there is no confusion, we will simply call it an orbifold bundle. We begin
with local uniformizing systems for orbifold bundles. Given a
uniformized topological space $U$ and a topological space $E$ with
a surjective continuous map $pr:E\rightarrow U$, a {\it
uniformizing system of rank $k$ orbifold bundle} for $E$ over $U$
consists of the following data:

\begin{itemize}
\item [{1.}] A uniformizing system $(V,G,\pi)$ of $U$.
\item [{2.}] A uniformizing system $(V\times\R^k,G,\tilde{\pi})$ for $E$.
The action of $G$ on $V\times\R^k$ is an extension of the action
of $G$ on $V$ given by $g\cdot (x,v)=(g\cdot x,\rho(x,g)v)$ where
$\rho:V\times G\rightarrow Aut(\R^k)$ is a smooth map satisfying:
$$ \rho(g\cdot x,h)\circ\rho(x,g)=\rho(x,hg), \hspace{3mm} g,h\in
G, x\in V. $$
\item [{3.}] The natural projection map $\tilde{pr}:V\times\R^k\rightarrow V$
satisfies $\pi\circ\tilde{pr}=pr\circ\tilde{\pi}$.
\end{itemize}

We can similarly define {\it isomorphisms} between uniformizing
systems of orbifold bundle for $E$ over $U$. The only additional
requirement is that the diffeomorphisms between $V\times\R^k$ are
linear on each fiber of $\tilde{pr}:V\times\R^k\rightarrow V$.
Moreover, for each connected open subset $U^\prime$ of $U$, we can
similarly prove that there is a unique isomorphism class of
induced uniformizing systems of orbifold bundle for
$E^\prime=pr^{-1}(U^\prime)$ over $U^\prime$. The {\it germ} of
uniformizing systems of orbifold bundle at a point $p\in U$ can be
also similarly defined.

\vspace{2mm}

\noindent{\bf Definition 2.3:}{\it
\begin{enumerate}
\item Let $X$ be an orbifold and $E$ be a topological space with a surjective
continuous map $pr:E\rightarrow X$. A {\it rank $k$ orbifold
bundle structure} on $E$ over $X$ consists of following data: For
each point $p\in X$, there is a uniformized neighborhood $U_p$ and
a uniformizing system of rank $k$ orbifold bundle for
$pr^{-1}(U_p)$ over $U_p$ such that for any $q\in U_p$, the
uniformizing systems of orbifold bundle over $U_p$ and $U_q$
define the same germ at $q$. The topological space $E$ with a given germ of
orbifold bundle
structures becomes an orbifold ($E$ is obviously Hausdorff and
second countable) and is called an {\it orbifold bundle} over $X$.
Each chart $(V_p\times\R^k,G_p,\tilde{\pi}_p)$ is called a {\it
local trivialization} of $E$. At each point $p\in X$, the fiber
$E_p=pr^{-1}(p)$ is isomorphic to $\R^k/G_p$. It contains a linear
subspace $E^p$ of fixed points of $G_p$.
\item The notion of {\it orbifold bundle over an orbifold with boundary} is
similarly defined. One can easily verify that if $pr:E\rightarrow
X$ is an orbifold bundle over an orbifold with boundary $X$, then
the restriction to the boundary $\partial X$, $E_{\partial
X}=pr^{-1}(\partial X)$, is an orbifold bundle over $\partial X$.
\item One can define fiber orbifold bundle in the same vein.
\end{enumerate}} \hfill $\Box$

A $C^l$ map $\tilde{s}$ from $X$ to an orbifold bundle
$pr:E\rightarrow X$ is called a {\it $C^l$ section} if locally
$\tilde{s}$ is given by $\tilde{s}_p: V_p\rightarrow
V_p\times\R^k$ where $\tilde{s}_p$ is $G_p$-equivariant and
$\tilde{pr}\circ\tilde{s}_p=Id$ on $V_p$. We observe that

\begin{itemize}
\item [{1.}] For each point $p$, $s(p)$ lies in $E^p$, the linear subspace
of fixed points of $G_p$.
\item [{2.}] The space of all $C^l$ sections of $E$, denoted by $C^l(E)$, has a
structure of vector space over $\R$ (or $\C$) as well as a
$C^l(X)$-module structure.
\item [{3.}] The $C^l$ sections $\tilde{s}$ are in $1:1$ correspondence with
the underlying continuous maps $s$.
\end{itemize}

Orbifold bundles are more conveniently described by transition maps, e.g. as in
\cite{S}. More precisely, an orbifold bundle over an orbifold $X$
can be constructed from the following data: A compatible cover
$\U$ of $X$ such that for any injection
$i:(V^\prime,G^\prime,\pi^\prime)\rightarrow (V,G,\pi)$, there is
a smooth map $g_i:V^\prime\rightarrow Aut(\R^k)$ giving an open
embedding $V^\prime\times\R^k\rightarrow V\times\R^k$ by
$(x,v)\rightarrow (i(x),g_i(x)v)$, and for any composition of
injections $j\circ i$, we have

$$ g_{j\circ i}(x)=g_j(i(x))\circ g_i(x).
\leqno(2.1) $$

Two collections of maps $g^{(1)}$ and $g^{(2)}$ define isomorphic
orbifold bundles if there are maps $\delta_V:V\rightarrow Aut(\R^k)$
such that for any injection $i: (V^\prime,G^\prime,\pi^\prime)
\rightarrow (V,G,\pi)$, we have
 $$g^{(2)}_i(x)=\delta_V(i(x))\circ g^{(1)}_i(x)\circ
(\delta_{V^\prime}(x))^{-1}, \forall x\in V^\prime.\leqno(2.2)$$

Since the equation (2.1) behaves naturally under constructions of
vector spaces such as tensor product, exterior product, etc., we
can define the corresponding constructions for orbifold bundles.

\vspace{2mm}

\noindent{\bf Example 2.4: }{\it  For an orbifold $X$, the tangent
bundle $TX$ can be constructed because the differential of any
injection satisfies the equation (2.1). Likewise, we define
cotangent bundle $T^\ast X$, the bundles of exterior power or
tensor product. The $C^\infty$ sections of these bundles give us
vector fields, differential forms or tensor fields on $X$. We
remark that if $\omega$ is a differential form on $X^\prime$ and
$\tilde{f}:X\rightarrow X^\prime$ is a $C^\infty$ map, then there
is a pull-back form $\tilde{f}^\ast\omega$ on $X$.} \vskip 0.1in

\vspace{2mm}

Let $U$ be an open subset of an orbifold $X$ with an orbifold
structure $\{(V_p,G_p,\pi_p):p\in X\}$, then
$\{(V^\prime_p,G^\prime_p,\pi^\prime_p): p\in U\}$ is an orbifold
structure on $U$ where $(V^\prime_p,G^\prime_p,\pi^\prime_p)$ is a
uniformizing system of $\pi_p(V_p)\cap U$ induced from
$(V_p,G_p,\pi_p)$. Likewise, let $pr:E\rightarrow X$ be an
orbifold bundle and $U$ an open subset of $X$, then
$pr:E_{U}=pr^{-1}(U)\rightarrow U$ inherits a unique germ of
orbifold bundle structures from $E$, called the {\it restriction
of $E$ over $U$}. When $U$ is a uniformized open set in $X$, say
uniformized by $(V,G,\pi)$, then there is a smooth vector bundle
$E_V$ over $V$ with a smooth action of $G$ such that
$(E_V,G,\tilde{\pi})$ uniformizes $E_U$. This is seen as follows:
We first take a compatible cover $\U$ of $U$, fine enough so that
the preimage under $\pi$ is a compatible cover of $V$. Let $E_U$
be given by a set of transition maps with respect to $\U$
satisfying (2.1), then the pull-backs under $\pi$ form a set of
transition maps with respect to $\pi^{-1}(\U)$ with an action of
$G$ by permutations, also satisfying (2.1), so that it defines a
smooth vector bundle over $V$ with a compatible smooth action of
$G$. Any $C^l$ section of $E$ on $X$ restricts to a $C^l$ section
of $E_U$ on $U$, and when $U$ is a uniformized open set by
$(V,G,\pi)$, it lifts to a $G$-equivariant $C^l$ section of $E_V$
on $V$.

\vspace{2mm}

Integration over orbifolds is defined as follows. Let $U$ be a
connected n-dimensional orbifold, which is uniformized by
$(V,G,\pi)$, with the kernel of the action of $G$ on $V$ denoted
by $K$. For any compact supported differential n-form $\omega$ on
$U$, which is, by definition, a $G$-equivariant compact supported
n-form $\tilde{\omega}$ on $V$, the integration of $\omega$ on $U$
is defined by
$$ \int_U^{orb} \omega:=\frac{1}{|G|}\int_V
\tilde{\omega}, \leqno (2.3)
$$
where $|G|$ is the order of the
group $G$. In general, let $X$ be an orbifold. Fix a $C^\infty$
partition of unity $\{\rho_i\}$ subordinated to $\{U_i\}$ where
each $U_i$ is a uniformized open set in $X$. Then the integration
over $X$ is defined by
$$ \int_X^{orb} \omega:=\sum_i
\int_{U_i}^{orb}\rho_i\;\omega, \leqno (2.4)
$$
which is independent of
the choice of the partition of unity $\{\rho_i\}$. We remark that it
is important throughout this paper that we adopt the integration over
orbifolds as in $(2.3)$ and $(2.4)$, where we divide the integral over
the uniformizer $V$ by the group order $|G|$ instead of $|G|/|K|$ ($K$
is the kernel of the action). As a result, the fundamental class of
an orbifold is rational in general. The integration $\int^{orb}$
coincides with the usual measure-theoretic integration if and only if
the orbifold is reduced.

\vspace{2mm}

The de Rham cohomology groups of an orbifold are defined similarly
through differential forms, which are naturally isomorphic to the
singular cohomology groups with real coefficients. For an
oriented, closed orbifold, the singular cohomology groups are
naturally isomorphic to the intersection homology groups, both
with rational coefficients, for which the Poincar\'{e} duality is
valid \cite{GM}.

Characteristic classes ({\it Euler class} for oriented orbifold
bundles, {\it Chern classes} for complex orbifold bundles, and
{\it Pontrjagin classes} for real orbifold bundles) are
well-defined for orbifold bundles. One way to define them is
through Chern-Weil theory, so that the characteristic classes take
values in the deRham cohomology groups. Another way to define them
is through the transgressions in the Serre spectral sequences with
rational coefficients of the associated Stiefel orbifold bundles,
so that these characteristic classes are defined over the
rationals \cite{K1}.

\section{Orbifold Cohomology Groups}

In this section, we introduce the main object of study, the {\it
orbifold cohomology groups} of an almost complex orbifold.

\subsection{Twisted sectors}

Let $X$ be an orbifold. For any point $p\in X$, let
$(V_p,G_p,\pi_p)$ be a local chart at $p$. Consider the set of
pairs:
$$
\widetilde{X}=\{(p,(g)_{G_p})|p\in X, g\in G_p\}, \leqno(3.1.1)
$$
where $(g)_{G_p}$ is the conjugacy class
of $g$ in $G_p$. If there is no confusion, we will omit the
subscript $G_p$ to simplify the notation. There is a surjective map
$\pi:\widetilde{X}\rightarrow X$ defined by $(p,(g))\mapsto p$.

\vspace{2mm}

\noindent{\bf Lemma 3.1.1 (Kawasaki,\cite{K1}): }{\it  The set
$\widetilde{X}$ is naturally an orbifold (not necessarily
connected) with an orbifold structure given by
$$
\{\pi_{p,g}:(V_p^{g},C(g))\rightarrow V_p^{g}/C(g): p\in X,
g\in G_p.\},
$$
where $V_p^{g}$ is the fixed-point set of $g$ in
$V_p$, $C(g)$ is the centralizer of $g$ in $G_p$. Moreover, if $X$
is closed, so is $\widetilde{X}$. Under this orbifold
structure, the map $\pi:\widetilde{X}\rightarrow X$ is a
$C^\infty$ map.}

\vspace{2mm}

\noindent{\bf Proof:} First we identify a point $(q,(h))$ in
$\widetilde{X}$ as a point in $\bigsqcup_{\{(g), g\in
G_p\}}V_p^g/C(g)$ if $q\in U_p$ for some $p\in X$.
Pick a representative $y\in V_p$ such that $\pi_p(y)=q$. Then this
gives rise to a monomorphism $\lambda_y:G_q\rightarrow G_p$. Pick
a representative $h\in G_q$ for $(h)$ in $G_q$, we let
$g=\lambda_y(h)$. Then $y\in V^g_p$. So we have a map
$\Phi:(q,h)\rightarrow (y,g)\in (V_p^g,G_p)$. If we change $h$ by
a $h^\prime =a^{-1}ha\in G_q$ for $a\in G_q$, then $g$ is changed
to $\lambda_y(a^{-1}ha) =\lambda_y(a)^{-1}g\lambda_y(a)$. So we have
$\Phi:(q,a^{-1}ha)\rightarrow (y,\lambda_y(a)^{-1}g\lambda_y(a))\in
(V^{\lambda_y(a)^{-1}g\lambda_y(a)}_p, G_p)$. (Note that $\lambda_y$
is determined up to conjugacy by an element in $G_q$.) If we take
a different representative $y^\prime\in V_p$ such that
$\pi_p(y^\prime)=q$, and suppose $y^\prime=b\cdot y$ for some
$b\in G_p$. Then we have a different identification
$\lambda_{y^\prime}:G_q\rightarrow G_p$ of $G_q$ as a subgroup of
$G_p$ where $\lambda_{y^\prime}=b\cdot\lambda_y\cdot b^{-1}$. In
this case, we have $\Phi:(q,h)\rightarrow (y^\prime,bgb^{-1})\in
(V_p^{bgb^{-1}},G_p)$. If $g=bgb^{-1}$, then $b\in C(g)$. In any
event, $\Phi$ induces a map $\phi$ sending $(q,(h))$ to a point in
$\bigsqcup_{\{(g), g\in G_p\}}V_p^g/C(g)$. It is one to one
because if $\phi(q_1,(h_1))=\phi(q_2,(h_2))$, then we may assume
that $\Phi(q_1,h_1)=\Phi(q_2,h_2)$ after applying conjugations.
But this means that $(q_1,h_1)=(q_2,h_2)$. It is easily seen that
this map $\phi$ is also onto. Hence we have shown that
$\widetilde{X}$ is covered by
$\bigsqcup_{\{p\in X\}}\bigsqcup_{\{(g), g\in G_p\}} V_p^g/C(g)$.

We define a topology on $\widetilde{X}$ so that each
$V_p^g/C(g)$ is an open subset for any $(p,g)$ where $p\in X$ and
$g\in G_p$. We also uniformize $V_p^g/C(g)$ by
$(V_p^g,C(g))$.  It remains to show that these charts fit together
to form an orbifold structure on $\widetilde{X}$. Let $x\in
V_p^g/C(g)$ and take a representative $\tilde{x}$ in $V_p^g$. Let
$H_x$ be the isotropy subgroup of $\tilde{x}$ in $C(g)$. Then
$(V_p^g,C(g))$ induces a germ of uniformizing system at $x$ as
$(B_x,H_x)$ where $B_x$ is a small ball in $V_p^g$ centered at
$\tilde{x}$. Let $\pi_p(\tilde{x})=q$. We need to write
$(B_x,H_x)$ as $(V_q^h,C(h))$ for some $h\in G_q$. We let
$\lambda_x: G_q\rightarrow G_p$ be an induced monomorphism
resulted from choosing $\tilde{x}$ as the representative of $q$ in
$V_p$. We define $h=\lambda_x^{-1}(g)$ ($g$ is in
$\lambda_x(G_q)$ since $\tilde{x} \in V^g_p$ and
$\pi_p(\tilde{x})=q$.) Then we can identify $B_x$ as $V_q^h$. We
also see that $H_x=\lambda_x(C(h))$. Therefore $(B_x,H_x)$ is
identified as $(V_q^h,C(h))$.

The map $\pi:\widetilde{X}\rightarrow X$ is obviously
continuous with the given topology of $\widetilde{X}$, and
actually is a $C^\infty$ map with the given orbifold structure on
$\widetilde{X}$ with the local liftings given by
embeddings $V_p^g\hookrightarrow V_p$.

We finish the proof by showing that $\widetilde{X}$ is
Hausdorff and second countable with the given topology. Let
$(p,(g))$ and $(q,(h))$ be distinct two points in
$\widetilde{X}$. When $p\neq q$, there are $U_p$, $U_q$
such that $U_p\cap U_q=\emptyset$ since $X$ is Hausdorff. It is
easily seen that in this case $(p,(g))$ and $(q,(h))$ are
separated by disjoint neighborhoods $\pi^{-1}(U_p)$ and
$\pi^{-1}(U_q)$, where $\pi:\widetilde{X}\rightarrow X$.
When $p=q$, we must then have $(g)\neq (h)$. In this case, $(p,(g))$ and
$(q,(h))$ lie in different open subsets $V_p^g/C(g)$ and
$V_q^h/C(h)$ respectively. Hence $\widetilde{X}$ is
Hausdorff. The second countability of $\widetilde{X}$
follows from the second countability of $X$ and the fact that
$\pi^{-1}(U_p)$ is a finite union of open subsets of
$\widetilde{X}$ for each $p\in X$ and a uniformized
neighborhood $U_p$ of $p$.
\hfill $\Box$

\vspace{2mm}

   Next, we would like
   to describe the connected components of $\widetilde{X}$. Recall
   that every point $p$ has a local chart $(V_p,G_p,\pi_p)$ which gives a
   local uniformized neighborhood $U_p=\pi_p(V_p)$.
   If $q\in U_p$, up to conjugation, there is an
   injective homomorphism $G_q\rightarrow G_p$. For $g\in G_q$,
   the conjugacy class $(g)_{G_p}$ is well-defined. We define
   an equivalence relation $(g)_{G_q}\sim (g)_{G_p}$. Let $T$ be
   the set of equivalence classes.  To abuse the notation, we
   often use $(g)$ to denote the equivalence class which $(g)_{G_q}$
   belongs to.
   It is clear that $\widetilde{X}$ is decomposed as a disjoint
   union of connected components
   $$
   \widetilde{X}=\bigsqcup_{(g)\in T} X_{(g)},\leqno(3.1.2)
   $$
   where
   $$
   X_{(g)}=\{(p,(g')_{G_p})|g'\in G_p, (g')_{G_p}\in (g)\}.\leqno(3.1.3)
   $$
   \vskip 0.1in

   \noindent{\bf Definition 3.1.2: }{\it $X_{(g)}$ for $g\neq 1$ is called
   a twisted sector.
   Furthermore, we call $X_{(1)}=X$ the nontwisted sector.}
   \vskip 0.1in

\noindent{\bf Example 3.1.3:} Consider the case that the orbifold
$X=Y/G$ is a global quotient. We will show that $\widetilde{X}$ can be
identified with $\bigsqcup_{\{(g),g\in G\}} Y^g/C(g)$ where $Y^g$
is the fixed-point set of element $g\in G$.

Let $\pi:\widetilde{X}\rightarrow X$ be the surjective map
defined by $(p,(g))\mapsto p$. Then for any $p\in X$, the
preimage $\pi^{-1}(p)$ in $\widetilde{X}$ has a neighborhood
described by $W_p=\bigsqcup_{\{(g),g\in G_p\}}V_p^g/C(g)$, which
is uniformized by $\widehat{W}_p=\bigsqcup_{\{(g),g\in G_p\}}
V_p^g$. For each $p\in X$, pick a $y\in Y$ that represents $p$,
and an injection $(\phi_p,\lambda_p): (V_p,G_p)\rightarrow (Y,G)$
whose image is centered at $y$. This induces an open embedding
$\tilde{f}_p:\widehat{W}_p\rightarrow
\bigsqcup_{\{(\lambda_p(g)),\lambda_p(g)\in G\}}Y^{\lambda_p(g)}
\subset \bigsqcup_{\{(g),g\in G\}}Y^g$, which induces a
homeomorphism $f_p$ from $W_p$ into $\bigsqcup_{\{(g),g\in
G\}}Y^g/C(g)$ that is independent of the choice of $y$ and
$(\phi_p,\lambda_p)$. These maps $\{f_p;p\in X\}$ fit together to
define a map $f:\widetilde{X}\rightarrow \bigsqcup_{\{(g),
g\in G\}}Y^g/C(g)$ which we can verify to be a homeomorphism.

\hfill $\Box$

\vspace{3mm}

\noindent {\bf Remark 3.1.4: }{\it There is a natural $C^\infty$
map $I:\widetilde{X}\rightarrow \widetilde{X}$
defined by $$ I((p,(g)_{G_p}))=(p, (g^{-1})_{G_p}).\leqno(3.1.4)
$$ The map $I$ is an involution (i.e., $I^2=Id$) which induces an
involution on the set $T$ of equivalence classes of relations
$(g)_{G_q}\sim  (g)_{G_p}$. We denoted by $(g^{-1})$ the image of
$(g)$ under this induced map.}

\vskip 0.1in

\subsection{Degree shifting and orbifold cohomology group}

For the rest of the paper, we will assume that $X$ is an almost
complex orbifold with an almost complex structure $J$. Recall that
an almost complex structure $J$ on $X$ is a smooth section of the
orbifold bundle $End(TX)$ such that $J^2=-Id$. Observe that
$\widetilde{X}$ naturally inherits an almost
complex structure from the one on $X$, and the map
$\pi:\widetilde{X}\rightarrow  X$ defined by
$(p,(g)_{G_p})\rightarrow p$ is naturally pseudo-holomorphic,
i.e., its differential commutes with the almost complex structures
on $\widetilde{X}$ and $X$.

An important feature of orbifold cohomology groups is degree
shifting, which we shall explain now. Let $p$ be any
point of $X$. The almost complex structure on $X$ gives rise to a
representation $\rho_p: G_p \rightarrow GL(n,\C)$ (here
$n=\dim_\C X$). For any $g\in G_p$, we write $\rho_p(g)$ as a
diagonal matrix
$$
diag(e^{2\pi i m_{1,g}/m_g}, \cdots, e^{2\pi i
m_{n,g}/m_g}),
$$
where $m_g$ is the order of $\rho_p(g)$, and
$0\leq m_{i,g} <m_g$. This matrix depends only on the conjugacy
class $(g)_{G_p}$ of $g$ in $G_p$. We define a function
$\iota:\widetilde{X}\rightarrow {\bf Q}$ by
$$
\iota(p,(g)_{G_p})=\sum_{i=1}^n \frac{m_{i,g}}{m_g}.
$$
It is straightforward to show the following

\vspace{2mm}

\noindent{\bf Lemma 3.2.1: }{\it The function $\iota:
X_{(g)}\rightarrow {\bf Q}$ is constant. Its constant value, which will
be denoted by $\iota_{(g)}$, satisfies the following
conditions:

\begin{itemize}
\item $\iota_{(g)}$ is integral if and only if $\rho_p(g)\in SL(n,\C)$.
\item
$$ \iota_{(g)}+\iota_{(g^{-1})}=rank(\rho_p(g)-I), \leqno(3.2.1)$$
which is the ``complex codimension'' $\dim_\C X-\dim_\C
X_{(g)}=n-\dim_\C X_{(g)}$ of $X_{(g)}$ in $X$. As a consequence,
$\iota_{(g)}+\dim_\C X_{(g)}<n$ when $\rho_p(g)\neq I$.
\end{itemize}}

\vskip 0.1in

\noindent{\bf Definition 3.2.2: }{\it $\iota_{(g)}$ is called a
degree shifting number.}

\vskip 0.1in

In the definition of orbifold cohomology groups, we will shift up
the degree of cohomology classes of $X_{(g)}$ by $2\iota_{(g)}$.
The reason for such a degree shifting will become clear after we
discuss the dimension of moduli space of ghost maps (see formula
(4.2.14)).

    An orbifold $X$ is called a {\em $SL$-orbifold} if
    $\rho_p(g)\in SL(n, \C)$ for all $p\in X$ and $g\in G_p$, and
    called a {\em $SP$-orbifold} if $\rho_p(g)\in SP(n,\C)$.
    In particular, a Calabi-Yau orbifold is a
    $SL$-orbifold, and a holomorphic symplectic orbifold or
    hyperkahler orbifold is a $SP$-orbifold. By Lemma 3.2.1,
    $\iota_{(g)}$ is integral if and only if $X$ is a
    $SL$-orbifold.

We observe that although the almost complex structure $J$ is
involved in the definition of degree shifting numbers
$\iota_{(g)}$, they do not depend on $J$ because locally the
parameter space of almost complex structures, which is the coset
$SO(2n,\R)/U(n,\C)$, is connected.

\vskip 0.1in

\noindent{\bf Definition 3.2.3: }{\it  We define the orbifold
cohomology groups $H^d_{orb}(X)$ of $X$  by
$$
H^d_{orb}(X)=\oplus_{(g)\in T}
H^{d-2\iota_{(g)}}(X_{(g)})\leqno(3.2.2)
$$
and orbifold Betti numbers $b^d_{orb}=\sum_{(g)}\dim H^{d-2\iota_{(g)}}
(X_{(g)})$.}

\vskip 0.1in

Here each $H^\ast(X_{(g)})$ is the singular cohomology of $X_{(g)}$
with real coefficients, which is isomorphic to the corresponding de Rham
cohomology group. As a consequence, the cohomology classes can be
represented by closed differential forms on $X_{(g)}$. Note that,
in general, orbifold cohomology groups are rationally graded.

\vspace{2mm}

Suppose $X$ is a  complex orbifold with an integrable complex
structure $J$. Then each twisted sector $X_{(g)}$ is also a
complex orbifold with the induced complex structure. We consider
the \v{C}ech cohomology groups on $X$ and each $X_{(g)}$ with
coefficients in the sheaves of holomorphic forms (in the orbifold
sense). These \v{C}ech cohomology groups are identified with the
Dolbeault cohomology groups of $(p,q)$-forms (in the orbifold
sense). When $X$ is closed,
the harmonic theory \cite{Ba} can be applied to show that
these groups are finite dimensional, and there is a Kodaira-Serre
duality between them. When $X$ is a closed Kahler orbifold (so is
each $X_{(g)}$), these groups are then related to the singular
cohomology groups of $X$ and $X_{(g)}$ as in the smooth case, and
the Hodge decomposition theorem holds for these cohomology groups.

\vspace{3mm}

\noindent{\bf Definition 3.2.4:} {\it Let $X$ be a  complex
orbifold. We define, for $0\leq p,q\leq \dim_{\C}X$, orbifold
Dolbeault cohomology groups
$$
H^{p,q}_{orb}(X)=\oplus_{(g)}
H^{p-\iota_{(g)},q-\iota_{(g)}} (X_{(g)}). \leqno(3.2.3)
$$
We define orbifold Hodge numbers by $h^{p,q}_{orb}(X)=\dim
H^{p,q}_{orb}(X)$.}

\vspace{2mm}
    \noindent
    {\bf Remark 3.2.5: }{\it We can define compact supported
    orbifold cohomology groups $H^*_{orb, c}(X),
    H^{*,*}_{orb, c}(X)$ in the obvious fashion.
    }
    \vskip 0.1in

\subsection{Poincar\'{e} duality}

Recall that there is a natural $C^\infty$ map
$I:X_{(g)}\rightarrow X_{(g^{-1})}$ defined by $(p,(g))\mapsto
(p,(g^{-1}))$, which is an automorphism of $\widetilde{X}$
as an orbifold and $I^2=Id$ (Remark 3.1.4).

\vspace{2mm}

\noindent{\bf Proposition 3.3.1: }(Poincar\'{e} duality){\it

For any $0\leq d\leq 2n$, the pairing
 $$
<\,\ >_{orb}:
H^d_{orb}(X)\times H^{2n-d}_{orb,c}(X)\rightarrow {\bf R}
$$
defined by the direct sum of
    $$
<\,\ >_{orb}^{(g)}:
H^{d-2\iota_{(g)}}(X_{(g)})\times H^{2n-d-2\iota_{(g^{-1})}}_c
(X_{(g^{-1})})\rightarrow {\bf R}
$$
where
    $$
<\alpha, \beta>^{(g)}_{orb} =\int_{X_{(g)}}^{orb} \alpha\wedge
I^*(\beta)\leqno(3.3.4)
$$
    for $\alpha\in H^{d-2\iota_{(g)}}(X_{(g)}), \beta\in
    H^{2n-d-2\iota_{(g^{-1})}}_c (X_{(g^{-1})})$
    is nondegenerate. Here the integral in the right hand side of $(3.3.4)$
    is defined using $(2.4)$.
}

\vskip 0.1in

Note that $<\,\ >_{orb}$ equals the ordinary Poincar\'{e} pairing
when restricted to the nontwisted sectors $H^\ast(X)$.

\vskip 0.1in

\noindent{\bf Proof:}
    By (3.2.1), we have
    $$
     2n-d-2\iota_{(g^{-1})}=\dim X_{(g)}-d-2\iota_{(g)}.
    $$
    Furthermore, $I|_{X_{(g)}}: X_{(g)}\rightarrow X_{(g^{-1})}$ is
    a homeomorphism.
    Under this homeomorphism, $<\,\ >^{(g)}_{orb}$
    is isomorphic to the ordinary Poincar\'{e} pairing on $X_{(g)}$.
    Hence $<\,\ >_{orb}$ is nondegenerate.
    \hfill $\Box$

\vspace{2mm}

For the case of orbifold Dolbeault cohomology, the following
proposition is straightforward.

\vspace{3mm}

\noindent{\bf Proposition 3.3.2:} {\it Let $X$ be an
$n$-dimensional  complex orbifold. There is a Kodaira-Serre
duality pairing $$ <\,\ >_{orb}:H^{p,q}_{orb}(X)\times
H^{n-p,n-q}_{orb,c}(X)\rightarrow \C $$ similarly defined as in
the previous proposition. When $X$ is closed and Kahler, the
following is true:
\begin{itemize}
\item $H^r_{orb}(X)\otimes\C=\oplus_{r=p+q}H^{p,q}_{orb}(X)$
\item $H^{p,q}_{orb}(X)=\overline{H^{q,p}_{orb}(X)}$,
\end{itemize}
and the two pairings (Poincar\'{e} and Kodaira-Serre) coincide.}

\section{Orbifold Cup Product and Orbifold Cohomology Ring}

\subsection{Orbifold cup product}
    In this section, we give an explicit definition of the orbifold cup
    product. Its interpretation in terms of Gromov-Witten invariants and the
proof of associativity of the product will be given in subsequent sections.

    Let $X$ be an orbifold, and $(V_p,G_p,\pi_p)$ be a uniformizing system
at point $p\in X$. We define the {\em $k$-multi-sector} of $X$, which is
denoted by $\widetilde{X}_k$, to be the set of all pairs $(p,(\g))$, where
$p\in X$, $\g=(g_1,\cdots,g_k)$ with each $g_i\in G_p$, and $(\g)$ stands
for the conjugacy class of $\g=(g_1,\cdots,g_k)$. Here two $k$-tuple
$(g_1^{(i)},\cdots,g_k^{(i)})$, $i=1,2$, are conjugate if there is a
$g\in G_p$
such that $g_j^{(2)}=gg_j^{(1)}g^{-1}$ for all $j=1,\cdots, k$.

\vspace{2mm}

\noindent{\bf Lemma 4.1.1:} {\it The $k$-multi-sector
$\widetilde{X}_k$ is naturally an orbifold, with the orbifold structure
given by
$$ \{\pi_{p,\g}:(V^\g_p,C(\g))\rightarrow V^\g_p/C(\g)\}, \leqno(4.1.1)
$$
where $V^\g_p=V_p^{g_1}\cap V_p^{g_2}\cap\cdots\cap V_p^{g_k}$,
$C(\g)=C(g_1)\cap C(g_2)\cap\cdots\cap C(g_k)$. Here
$\g=(g_1,\cdots,g_k)$, $V_p^g$ stands for the fixed-point set of
$g\in G_p$ in $V_p$, and $C(g)$ for the centralizer of $g$ in
$G_p$. For each $i=1,\cdots, k$, there is a $C^\infty$ map
$e_i:\widetilde{X}_k \rightarrow \widetilde{X}$ defined by sending
$(p,(\g))$ to $(p,(g_i))$ where $\g=(g_1,\cdots,g_k)$.
When $X$ is almost complex, $\widetilde{X}_k$ inherits an almost complex
structure from $X$, and when $X$ is closed, $\widetilde{X}_k$ is a finite
disjoint union of closed orbifolds.}

\vspace{2mm}

\noindent{\bf Proof:} The proof is parallel to the proof of Lemma
3.1.1 where $\widetilde{X}$ is shown to be an orbifold.

First we identify a point $(q,(\h))$ in $\widetilde{X}_k$
as a point in $\bigsqcup_{\{(p,(\g))\in \widetilde{X}_k\}}
V_p^{\g}/C(\g)$ if $q\in U_p$. Pick a representative $y\in
V_p$ such that $\pi_p(y)=q$. Then this gives rise to a
monomorphism $\lambda_y:G_q\rightarrow G_p$. Pick a representative
$\h=(h_1,\cdots,h_k)\in G_q\times \cdots\times G_q$ for $(\h)$, we
let $\g=\lambda_y(\h)$. Then $y\in V^{\g}_p$. So we have a map
$\theta:(q,\h)\rightarrow (y,\g)$. If we change $\h$ by
$\h^\prime=a^{-1}\h a$ for some $a\in G_q$, then $\g$ is changed
to $\lambda_y(a^{-1}\h a)=\lambda_y(a)^{-1}\g\lambda_y(a)$. So we have
$\theta:(q,a^{-1}\h a)\rightarrow (y,\lambda_y(a)^{-1}\g\lambda_y(a))$
where $y$ is regarded as a point in $V^{\lambda_y(a)^{-1}\g\lambda_y(a)}_p$.
(Note that $\lambda_y$ is determined up to conjugacy by an element in $G_q$.)
If we take a different representative $y^\prime\in V_p$ such that
$\pi_p(y^\prime)=q$, and suppose $y^\prime=b\cdot y$ for some
$b\in G_p$. Then we have a different identification
$\lambda_{y^\prime}:G_q\rightarrow G_p$ of $G_q$ as a subgroup of
$G_p$ where $\lambda_{y^\prime}=b\cdot\lambda_y\cdot b^{-1}$. In
this case, we have $\theta:(q,\h)\rightarrow (y^\prime,b\g
b^{-1})$ where $y^\prime\in V_p^{b\g b^{-1}}$. If $\g=b\g b^{-1}$,
then $b\in C(\g)$. Therefore we have shown that $\theta$ induces a
map sending $(q,(\h))$ to a point in
$\bigsqcup_{\{(p,(\g))\in\widetilde{X}_k\}}V_p^{\g}/C(\g)$,
which can be similarly shown to be one to one and onto. Hence we
have shown that $\widetilde{X}_k$ is covered by
$\bigsqcup_{\{(p,(\g))\in\widetilde{X}_k\}}V_p^{\g}/C(\g)$.

We define a topology on $\widetilde{X}_k$ so that each
$V_p^{\g}/C(\g)$ is an open subset for any $(p,\g)$. We also
uniformize $V_p^{\g}/C(\g)$ by $(V_p^{\g},C(\g))$. It remains to
show that these charts fit together to form an orbifold structure
on $\widetilde{X}_k$. Let $x\in V_p^{\g}/C(\g)$ and take a
representative $\tilde{x}$ in $V_p^\g$. Let $H_x$ be the isotropy
subgroup of $\tilde{x}$ in $C(\g)$. Then $(V_p^\g,C(\g))$ induces
a germ of uniformizing system at $x$ as $(B_x,H_x)$ where $B_x$ is
a small ball in $V_p^\g$ centered at $\tilde{x}$. Let
$\pi_p(\tilde{x})=q$. We need to write $(B_x,H_x)$ as
$(V_q^\h,C(\h))$ for some $\h\in G_q\times\cdots\times G_q$.
We let $\lambda_x: G_q\rightarrow G_p$ be an induced monomorphism
resulted from choosing $\tilde{x}$ as the representative of $q$ in
$V_p$. We define $\h=\lambda_x^{-1}(\g)$ (each $g_i$ is in
$\lambda_x(G_q)$ since $\tilde{x}\in V^\g_p$ and
$\pi_p(\tilde{x})=q$.) Then we can identify $B_x$ as $V_q^\h$. We
also see that $H_x=\lambda_x(C(\h))$. Therefore $(B_x,H_x)$ is
identified as $(V_q^\h,C(\h))$. Hence we proved that
$\widetilde{X}_k$ is naturally an orbifold with the orbifold
structure described above ($\widetilde{X}_k$ is Hausdorff
and second countable with the given topology for similar reasons).
The rest of the lemma is obvious.
\hfill $\Box$

\vspace{2mm}
    We can also describe the components of $\widetilde{X}_k$ in the same
    fashion. Using the conjugacy class of monomorphisms
    $\pi_{pq}: G_q\rightarrow G_p$ in the patching condition, we can define
    an equivalence relation $(\g)_{G_q}\sim (\pi_{pq}(\g))_{G_p}$ similarly.
    Let $T_k$ be the set of equivalence classes. We will write a general
    element of $T_k$ as $(\g)$. Then $\widetilde{X}_k$ is decomposed as a
disjoint union of connected orbifolds
$$
\widetilde{X}_k=\bigsqcup_{(\g)\in T_k} X_{(\g)},\leqno(4.1.2)
$$
where
$$
X_{(\g)}=\{(p, (\g')_{G_p})| (\g')_{G_p}\in (\g)\}.\leqno(4.1.3)
$$
There is a map $o:T_k\rightarrow T$ induced by the map
$o:(g_1,g_2,\cdots,g_k)\mapsto g_1g_2\cdots g_k$.
We set $T_k^o=o^{-1}((1))$. Then $T_k^o\subset T_k$ is the subset
of equivalence classes $(\g)$ such that $\g=(g_1,\cdots,g_k)$
satisfies the condition $g_1\cdots g_k=1$. Finally, we set
$$
\widetilde{X}_k^o:=\bigsqcup_{(\g)\in T_k^o} X_{(\g)}.\leqno(4.1.4)
$$

In order to define the orbifold cup product, we need a digression on a few
classical results about {\it reduced} 2-dimensional orbifolds
(cf. \cite{Th}, \cite{Sc}). Every closed orbifold of dimension 2 is complex,
whose underlying topological space is a closed Riemann surface.
More concretely, a closed, reduced 2-dimensional orbifold consists of the
following data: a closed Riemann surface $\Sigma$ with complex structure $j$,
a finite subset of distinct points $\z= (z_1,\cdots,z_k)$ on
$\Sigma$, each with a multiplicity $m_i\geq 2$ (let
$\m=(m_1,\cdots,m_k)$), such that the orbifold structure at $z_i$
is given by the ramified covering $z\rightarrow z^{m_i}$. We will
also call a closed, reduced 2-dimensional orbifold a {\it complex
orbicurve} when the underlying complex analytic structure is emphasized.

A $C^\infty$ map $\tilde{\pi}$ between two reduced connected
2-dimensional orbifolds is called an {\it orbifold covering} if the local
liftings of $\tilde{\pi}$
are either a diffeomorphism or a ramified covering. It is shown that the
universal orbifold covering exists, and its group of deck transformations
is defined to be the {\it orbifold fundamental group} of the orbifold. In fact,
given a reduced 2-orbifold $\Sigma$, with orbifold fundamental group denoted
by $\pi_1^{orb}(\Sigma)$, for any subgroup $\Gamma$ of $\pi_1^{orb}(\Sigma)$,
there is a reduced 2-orbifold $\widetilde{\Sigma}$ and an orbifold covering
$\tilde{\pi}:\widetilde{\Sigma}\rightarrow\Sigma$ such that $\tilde{\pi}$
induces an injective homomorphism $\pi_1^{orb}(\widetilde{\Sigma})\rightarrow
\pi_1^{orb}(\Sigma)$ with image $\Gamma\subset\pi_1^{orb}(\Sigma)$. The
orbifold fundamental group of a reduced, closed 2-orbifold $(\Sigma,\z,\m)$
has a presentation
$$
\pi_1^{orb}(\Sigma)=\{x_i,y_i,\lambda_j,i=1,\cdots,g,j=1,\cdots,k|
\prod_ix_iy_ix_i^{-1}y_i^{-1}\prod_j\lambda_j=1, \lambda_j^{m_j}=1\},
$$
where $g$ is the genus of $\Sigma$, $\z=(z_1,\cdots,z_k)$ and
$\m=(m_1,\cdots,m_k)$.

\vspace{2mm}

The remaining ingredient is to construct an ``obstruction bundle'' $E_{(\g)}$
over each component $X_{(\g)}$ where $(\g)\in T_3^o$. For this purpose,
we consider the Riemann sphere $S^2$ with three distinct marked points
$\z=(0,1,\infty)$. Suppose $(\g)$ is represented by $\g=(g_1,g_2,g_3)$ and
the order of $g_i$ is $m_i$ for $i=1,2,3$. We give a reduced orbifold
structure on $S^2$ by assigning $\m=(m_1,m_2,m_3)$ as the multiplicity of $\z$.
The orbifold fundamental group $\pi_1^{orb}(S^2)$ has the following
presentation
$$
\pi^{orb}_1(S^2)=\{\lambda_1, \lambda_2, \lambda_3|
   \lambda^{m_i}_i=1, \lambda_1\lambda_2\lambda_3=1\},
$$
where each generator $\lambda_i$ is geometrically represented by a loop
around the marked point $z_i$ (here recall that $(z_1,z_2,z_3)=(0,1,\infty)$).

Now for each point $(p,(\g)_{G_p})\in X_{(\g)}$, fix a representation
$\g$ of $(\g)_{G_p}$ where $\g=(g_1,g_2,g_3)$, we define a homomorphism
$\rho_{p,\g}: \pi^{orb}_1(S^2)\rightarrow G_p$ by sending $\lambda_i$ to
$g_i$, which is possible since $g_1g_2g_3=1$. Let $G\subset G_p$ be the image
of $\rho_{p,\g}$. There is a reduced 2-orbifold $\Sigma$ and an orbifold
covering $\tilde{\pi}:\Sigma\rightarrow S^2$, which induces the following
short exact sequence
$$
1\rightarrow\pi_1(\Sigma)\rightarrow\pi_1^{orb}(S^2)\rightarrow G\rightarrow 1.
$$
The group $G$ acts on $\Sigma$ as the group of deck transformations, whose
finiteness implies that $\Sigma$ is closed. Moreover, $\Sigma$ actually has
a trivial orbifold structure (i.e. $\Sigma$ is a Riemann surface)
since each map $\lambda_i\mapsto g_i$ is injective, and we can assume $G$ acts
on $\Sigma$ holomorphically. At end, we obtained a uniformizing system
$(\Sigma,G,\tilde{\pi})$ of $(S^2,\z,\m)$, which depends on $(p,\g)$, but
is locally constant.

The ``obstruction bundle'' $E_{(\g)}$ over $X_{(\g)}$ is constructed as
follows. On the local chart $(V_p^\g,C(\g))$ of $X_{(g)}$, $E_{(\g)}$ is
given by $(H^1(\Sigma)\otimes TV_p)^G\times V_p^\g\rightarrow V_p^\g$, where
$(H^1(\Sigma)\otimes TV_p)^G$ is the invariant subspace of $G$. We define an
action of $C(\g)$ on $H^1(\Sigma)\otimes TV_p$, which is trivial on the first
factor and the usual one on $TV_p$, then it is clear that $C(\g)$ commutes
with $G$, hence $(H^1(\Sigma)\otimes TV_p)^G$ is invariant under $C(\g)$.
In summary, we have obtained an action of $C(\g)$ on
$(H^1(\Sigma)\otimes TV_p)^G\times V_p^\g\rightarrow V_p^\g$, extending the
usual one on $V_p^\g$, and it is easily seen that these trivializations fit
together to define the bundle $E_{(\g)}$ over $X_{(\g)}$. If we set $e:
X_{(\g)}\rightarrow X$ to be the $C^\infty$ map $(p,(\g)_{G_p})\mapsto p$,
one may think of $E_{(\g)}$ as $(H^1(\Sigma)\otimes e^\ast TX)^G$.

\vspace{1.5mm}

Since we do not assume that $X$ is compact, $X_{(\g)}$ could be a non-compact
orbifold in general. The Euler class of $E_{(\g)}$ depends on a choice of
connection on $E_{(\g)}$.  Let $e_A(E_{(\g)})$ be the Euler form computed
from connection $A$ by Chern-Weil theory.

\vskip 0.1in

\noindent{\bf Definition 4.1.2: }{\it For $\alpha,\beta \in H^\ast_{orb}(X)$,
and $\gamma\in H^\ast_{orb,c}(X)$, we define a 3-point function
$$
<\alpha, \beta, \gamma>_{orb}=\sum_{(\g)\in T^0_3}\int_{X_{(\g)}}^{orb}
    e_1^\ast\alpha\wedge e_2^\ast\beta\wedge e^\ast_3\gamma \wedge
    e_A(E_{(\g)}),\leqno(4.1.5)
$$
where each $e_i:X_{(\g)}\rightarrow \widetilde{X}$ is the $C^\infty$ map
defined by $(p,(\g)_{G_p})\mapsto (p,(g_i)_{G_p})$ for $\g=(g_1,g_2,g_3)$.
Integration over orbifolds is defined by equation $(2.4)$.
}
\vskip 0.1in

    Note that since $\gamma$ is compact supported, each
    integral is finite, and the summation is over a finite subset of
    $T_3^o$. Moreover, if we choose different connection $A'$,
    $e_A(E_{(\g)}), e_{A'}(E_{(\g)})$ differ by an exact form. Hence the
    3-point function is independent of the choice of the connection $A$.

    \vskip 0.1in
    \noindent
    {\bf Definition 4.1.3: }{ \it We define the orbifold cup product
    on $H^\ast_{orb}(X)$ by the relation
    $$
    <\alpha\cup_{orb} \beta, \gamma>_{orb}=<\alpha, \beta,
    \gamma>_{orb}.\leqno(4.1.6)
    $$
    }

Next we shall give a decomposition of the orbifold cup product
$\alpha\cup_{orb} \beta$ according to the decomposition $H^\ast_{orb}(X)=
\oplus_{(g)\in T}H^{\ast-2\iota_{(g)}}(X_{(g)})$, when $\alpha,\beta$ are
homogeneous, i.e. $\alpha\in H^\ast(X_{(g_1)})$ and $\beta\in
H^\ast(X_{(g_2)})$ for some $(g_1),(g_2)\in T$. We need to introduce some
notation first. Given $(g_1),(g_2)\in T$, let $T((g_1),(g_2))$ be the subset
of $T_2$ which consists of $(\h)$ where $\h=(h_1,h_2)$ satisfies $(h_1)=(g_1)$
and $(h_2)=(g_2)$. Recall that there is map $o:T_k\rightarrow T$ defined by
sending $(g_1,g_2,\cdots,g_k)$ to $g_1g_2\cdots g_k$. We define a map
$\delta: \g\mapsto (\g,o(\g)^{-1})$, which clearly induces a one to one
correspondence between $T_k$ and $T_{k+1}^o$. We also denote by $\delta$ the
resulting isomorphism $\widetilde{X}_k\cong \widetilde{X}_{k+1}^o$. Finally,
we set $\delta_i=e_i\circ\delta$.

\vskip 0.1in

\noindent
    {\bf Decomposition Lemma 4.1.4: }{\it For any
$\alpha\in H^\ast(X_{(g_1)})$, $\beta\in H^\ast(X_{(g_2)})$,
    $$
\alpha\cup_{orb}\beta =\sum_{(\h)\in T((g_1),(g_2))}
(\alpha\cup_{orb}\beta)_{(\h)},\leqno(4.1.7)
$$
    where $(\alpha\cup_{orb}\beta)_{(\h)}\in
    H^*(X_{o((\h))})$ is defined by the relation
    $$
<(\alpha\cup_{orb}\beta)_{o((\h))}, \gamma>_{orb}
    =\int_{X_{(\h)}}^{orb}\delta_1^\ast\alpha\wedge \delta_2^\ast\beta\wedge
\delta^\ast_3\gamma\wedge e_A(\delta^\ast E_{\delta(\h)}),\leqno(4.1.8)
$$
for $\gamma\in H^*_c(X_{(o(\h)^{-1})})$.
}

\vspace{2mm}

In the subsequent sections, we shall describe the 3-point function and
orbifold cup product in terms of Gromov-Witten invariants. In fact,
we will prove the following

\vspace{2mm}

\noindent{\bf Theorem 4.1.5: } {\it Let $X$ be an almost complex
orbifold with almost complex structure $J$ and $\dim_\C X=n$. The orbifold
cup product preserves the orbifold grading, i.e.,
$$
\cup_{orb}: H_{orb}^p(X)\times H_{orb}^q(X) \rightarrow
H_{orb}^{p+q}(X)
$$
for any $0\leq p,q\leq 2n$ such that $p+q\leq
2n$, and has the following properties:
\begin{enumerate}
\item The total orbifold cohomology group $H^\ast_{orb}(X)=
\oplus_{0\leq d\leq 2n}H^d_{orb}(X)$ is a ring with unit
$e_X^0\in H^0(X)$ under $\cup_{orb}$, where $e_X^0$ is the
Poincar\'{e} dual to the fundamental class $[X]$. In particular,
$\cup_{orb}$ is associative.
\item When $X$ is closed, for each $H_{orb}^d(X)\times H_{orb}^{2n-d}(X)
\rightarrow H_{orb}^{2n}(X)$, we have
$$
\int_X^{orb} \alpha\cup_{orb}
\beta=<\alpha, \beta>_{orb}. \leqno(4.1.9)
$$
\item The cup product $\cup_{orb}$ is invariant under deformation of
$J$.
\item When $X$ is of integral degree shifting numbers, the total orbifold
cohomology group $H_{orb}^\ast(X)$ is integrally graded, and we
have supercommutativity $$
\alpha_1\cup_{orb}\alpha_2=(-1)^{\deg\alpha_1\cdot\deg\alpha_2}\alpha_2
\cup_{orb}\alpha_1. $$
\item Restricted to the nontwisted sectors, i.e., the ordinary
cohomologies $H^\ast(X)$, the cup product $\cup_{orb}$ equals
the ordinary cup product on $X$.
\end{enumerate}
}

\vspace{2mm}

When $X$ is a complex orbifold, the definition of orbifold cup product
$\cup_{orb}$ on the total orbifold Dolbeault cohomology group of $X$ is
completely parallel. We observe that in this case all the objects we have
been dealing with are holomorphic, i.e., $\widetilde{X}_k$ is a complex
orbifold, the ``obstruction bundles'' $E_{(\g)}\rightarrow X_{(\g)}$ are
holomorphic orbifold bundles, and the evaluation maps $e_i$ are holomorphic.

\vspace{2mm}

\noindent{\bf Definition 4.1.6: }{\it For any $\alpha_1\in
H^{p,q}_{orb}(X)$, $\alpha_2\in H_{orb}^{p^\prime,q^\prime}(X)$,
we define a 3-point function and orbifold cup product in the same fashion
as in Definitions 4.1.2, 4.1.3.}
\hfill $\Box$

\vskip 0.1in

Note that since the top Chern class of a holomorphic orbifold
bundle can be represented by a closed $(r,r)$-form where $r$ is
the (complex) rank of the bundle, it follows that the orbifold cup product
preserves the orbifold bi-grading, i.e.,
$\cup_{orb}:H_{orb}^{p,q}(X)\times
H_{orb}^{p^\prime,q^\prime}(X)\rightarrow H_{orb}^{p+p^\prime,q+q^\prime}(X)$.

\vspace{2mm}

The following theorem can be similarly proved.

\vspace{2mm}

\noindent{\bf Theorem 4.1.7: }{\it Let $X$ be a n-dimensional
complex orbifold with complex structure $J$. The orbifold
cup product
$$
\cup_{orb}: H_{orb}^{p,q}(X)\times
H_{orb}^{p^\prime,q^\prime}(X) \rightarrow
H_{orb}^{p+p^\prime,q+q^\prime}(X)
$$
has the following properties:
\begin{enumerate}
\item The total orbifold Dolbeault cohomology group is a ring with unit
$e_X^0\in H_{orb}^{0,0}(X)$ under $\cup_{orb}$, where $e_X^0$
is the class represented by the equaling-one constant function on
$X$.
\item When $X$ is closed, for each $H_{orb}^{p,q}(X)\times H_{orb}^{n-p,n-q}(X)
\rightarrow H_{orb}^{n,n}(X)$, the integral
$\int_{X}\alpha\cup_{orb} \beta$ equals the Kodaira-Serre pairing
$<\alpha,\beta>_{orb}$.
\item The cup product $\cup_{orb}$ is invariant under
deformation of $J$.
\item When $X$ is of integral degree shifting numbers, the total orbifold
Dolbeault cohomology group of $X$ is integrally graded, and we
have supercommutativity $$
\alpha_1\cup_{orb}\alpha_2=(-1)^{\deg\alpha_1\cdot\deg\alpha_2}\alpha_2
\cup_{orb}\alpha_1. $$
\item Restricted to the nontwisted sectors, i.e., the ordinary Dolbeault
cohomologies $H^{\ast,\ast}(X)$, the cup product $\cup_{orb}$
coincides with the ordinary wedge product on $X$.
\item When $X$ is Kahler and closed, the cup product $\cup_{orb}$ coincides
with the orbifold cup product on the total orbifold cohomology group
$H^\ast_{orb}(X)$ under the relation
$$
H^r_{orb}(X)\otimes\C=\oplus_{p+q=r}H^{p,q}_{orb}(X).
$$
\end{enumerate}
}

\subsection{Moduli space of ghost maps}

We first give a classification of rank-n complex orbifold bundles
over a closed, {\it reduced}, 2-dimensional orbifold.

Let $(\Sigma,\z,\m)$ be a closed, reduced, 2-dimensional orbifold,
where $\z=(z_1,\cdots,z_k)$ and $\m=(m_1,\cdots,m_k)$. Let $E$ be a
complex orbifold bundle of rank $n$ over $(\Sigma,\z,\m)$. Then at
each singular point $z_i$, $i=1,\cdots,k$, $E$ determines a
representation $\rho_i:\Z_{m_i}\rightarrow Aut(\C^n)$ so that over
a disc neighborhood of $z_i$, $E$ is uniformized by
$(D\times\C^n,\Z_{m_i},\pi)$ where the action of $\Z_{m_i}$ on
$D\times\C^n$ is given by $$ e^{2\pi i/m_i}\cdot (z,w)=(e^{2\pi
i/m_i}z,\rho_i(e^{2\pi i/m_i})w) \leqno(4.2.1)$$ for any
$w\in\C^n$. Each representation $\rho_i$ is uniquely determined by
a $n$-tuple of integers $(m_{i,1},\cdots,m_{i,n})$ with $0\leq
m_{i,j}<m_i$, as it is given by matrix $$ \rho_i(e^{2\pi i/m_i})=
diag(e^{2\pi im_{i,1}/m_i},\cdots, e^{2\pi im_{i,n}/m_i}).
\leqno(4.2.2)$$ Over the punctured disc $D_i\setminus\{0\}$ at
$z_i$, $E$ inherits a specific trivialization from
$(D\times\C^n,\Z_{m_i},\pi)$ as follows: We define a
$\Z_{m_i}$-equivariant map
$\Psi_i:D\setminus\{0\}\times\C^n\rightarrow
D\setminus\{0\}\times\C^n$ by $$ (z,w_1,\cdots,w_n)\rightarrow
(z^{m_i}, z^{-m_{i,1}}w_1,\cdots,z^{-m_{i,n}}w_n), \leqno(4.2.3)$$
where $Z_{m_i}$ acts trivially on the second
$D\setminus\{0\}\times\C^n$. Hence $\Psi_i$ induces a
trivialization $\psi_i: E_{D_i\setminus\{0\}} \rightarrow
D_i\setminus\{0\}\times\C^n$. We can extend the smooth complex
vector bundle $E_{\Sigma\setminus\z}$ over $\Sigma\setminus\z$ to
a smooth complex vector bundle over $\Sigma$ by using these
trivializations $\psi_i$. We call the resulting complex vector
bundle the {\it de-singularization} of $E$, and denote it by
$|E|$.

\vspace{2mm}

\noindent{\bf Proposition 4.2.1: }{\it The space of isomorphism
classes of complex orbifold bundles of rank $n$ over a closed, reduced,
2-dimensional orbifold $(\Sigma,\z,\m)$ where
$\z=(z_1,\cdots,z_k)$ and $\m=(m_1,\cdots,m_k)$, is in 1:1
correspondence with the set of
$(c,(m_{1,1},\cdots,m_{1,n}),\cdots,(m_{k,1},\cdots,m_{k,n}))$ for
$c\in{\bf Q}$, $m_{i,j}\in\Z$, where $c$ and $m_{i,j}$ are confined by
the following condition: $$ 0\leq m_{i,j}<m_i \hspace{2mm}\mbox{
and }  \hspace{2mm} c\equiv
\sum_{i=1}^{k}\sum_{j=1}^{n}\frac{m_{i,j}}{m_i} \pmod{\Z}.
\leqno(4.2.4)$$ In fact, $c$ is the first Chern number of the
orbifold bundle and
$c-(\sum_{i=1}^{k}\sum_{j=1}^{n}\frac{m_{i,j}}{m_i})$ is the first
Chern number of its de-singularization.}

\vspace{2mm}

\noindent{\bf Proof:} We only need to show the relation: $$
c_1(E)([\Sigma])=c_1(|E|)([\Sigma])+
\sum_{i=1}^{k}\sum_{j=1}^{n}\frac{m_{i,j}}{m_i}. \leqno(4.2.5)$$
We take a connection $\nabla_0$ on $|E|$ which equals $d$ on a
disc neighborhood $D_i$ of each $z_i\in \z$ so that
$c_1(|E|)([\Sigma])= \int_{\Sigma}c_1(\nabla_0)$. We use
$\nabla_0^\prime$ to denote the pull-back connection
$br_i^\ast\nabla_0$ on $D\setminus\{0\}\times \C^n$ via
$br_i:D\rightarrow D_i$ by $z\rightarrow z^{m_i}$. On the other
hand, on each uniformizing system $(D\times\C^n,\Z_{m_i},\pi)$, we
take the trivial connection $\nabla_i=d$ which is obvious
$\Z_{m_i}$-equivariant. Furthermore, we take a
$\Z_{m_i}$-equivariant cut-off function $\beta_i$ on $D$ which
equals one in a neighborhood of the boundary $\partial D$. We are
going to paste these connections together to get a connection
$\nabla$ on $E$. We define $\nabla$ on each uniformizing system
$(D\times\C^n,\Z_{m_i},\pi)$ by $$ \nabla_v
u=(1-\beta_i)(\nabla_i)_v u+\beta_i\bar{\psi}_i^{-1}(\nabla_0)
_{\bar{\psi}_i v}\bar{\psi}_i u, \leqno(4.2.6)$$ where
$\bar{\psi}_i: D\setminus\{0\}\times\C^n\rightarrow
D\setminus\{0\}\times \C^n$ is given by $$
(z,w_1,\cdots,w_n)\rightarrow (z,
z^{-m_{i,1}}w_1,\cdots,z^{-m_{i,n}}w_n). \leqno(4.2.7)$$ One
easily verifies that $F(\nabla)=F(\nabla_0)$ on $\Sigma\setminus
(\cup_i D_i)$ and $$F(\nabla)=-diag(d(\beta_i m_{i,1}dz/z),\cdots,
d(\beta_i m_{i,n}dz/z))$$ on each uniformizing system
$(D,\Z_{m_i},\pi)$. So
\begin{eqnarray*}
c_1(E)([\Sigma]) & = & \int_{\Sigma}^{orb}c_1(\nabla)\\
                 & = & \int_{\Sigma\setminus(\cup_i D_i)}c_1(\nabla_0)+
                       \sum_{i=1}^k\frac{1}{m_i}\int_{D}c_1(\nabla)\\
                 & = & c_1(|E|)([\Sigma])+
                       \sum_{i=1}^{k}\sum_{j=1}^{n}\frac{m_{i,j}}{m_i}.
\end{eqnarray*}
Here the integraton over $\Sigma$, $\int_{\Sigma}^{orb}$, should be
understood as in $(2.4)$.
\hfill $\Box$

We will need the following index formula.

\vspace{2mm}

\noindent
    {\bf Proposition 4.2.2: }
    {\it Let $E$ be a holomorphic orbifold bundle of rank $n$
    over a complex orbicurve $(\Sigma,\z,\m)$ of genus $g$. Then
    ${\cal O}(E)={\cal O}(|E|)$, where ${\cal O}(E), {\cal O}(|E|)$
    are sheaves of holomorphic sections of $E, |E|$. Hence,
    $$ \chi({\cal O}(E))=\chi({\cal O}(|E|))=c_1(|E|)([\Sigma])+n(1-g).
    \leqno(4.2.9) $$
    If $E$ corresponds to $(c,(m_{1,1},\cdots,m_{1,n}),\cdots,
    (m_{k,1},\cdots,m_{k,n}))$ (cf. Proposition 4.2.1), then we have
    $$\chi({\cal O}(E))=n(1-g)+c_1(E)([\Sigma])-\sum_{i=1}^k\sum_{j=1}^n
      \frac{m_{i,j}}{m_i}.$$ }
    \vskip 0.1in

\vspace{2mm}

\noindent{\bf Proof:} By construction, we have  ${\cal O}(E)={\cal
O}(|E|)$. Hence $$\chi({\cal O}(E))=\chi({\cal
O}(|E|))=c_1(|E|)([\Sigma])+n(1-g).
    \leqno(4.2.10) $$ By proposition 4.2.1, we have
$$\chi({\cal
O}(E))=n(1-g)+c_1(E)([\Sigma])-\sum_{i=1}^k\sum_{j=1}^n
      \frac{m_{i,j}}{m_i},$$
if  $E$ corresponds to $(c,(m_{1,1},\cdots,m_{1,n}),\cdots,
    (m_{k,1},\cdots,m_{k,n}))$. \hfill $\Box$

\vspace{3mm}

Now we come to the main issue of this section. Recall that suppose
$f:X \rightarrow X^\prime$ is a $C^\infty$ map between manifolds
and $E$ is a smooth vector bundle over $X^\prime$, then there is a
smooth pull-back vector bundle $f^\ast E$ over $X$ and a bundle
morphism $\bar{f}:f^\ast E\rightarrow E$ which covers the map $f$.
However, if instead,  we have a $C^\infty$ map $\tilde{f}$ between
orbifolds $X$ and $X^\prime$, and an orbifold bundle $E$ over
orbifold $X^\prime$, the question whether there is a pull-back
orbifold bundle $E^\ast$ over $X^\prime$ and an orbifold bundle
morphism $\bar{f}:E^\ast\rightarrow E$ covering the map
$\tilde{f}$ is a quite complicated issue: (1) What is the precise meaning of
pull-back orbifold bundle $E^\ast$, (2) $E^\ast$ might not
exist, or even if it exists, it might not be unique. Understanding this
question is the first step in our establishment of an {\it orbifold
Gromov-Witten theory} in \cite{CR}.

In the present case, given a constant map $f:\Sigma\rightarrow X$
from a marked Riemann surface $\Sigma$ with marked-point set $\z$
into an almost complex orbifold $X$, we need to settle the
existence and classification problem of pull-back orbifold bundles
via $f$, with some reduced orbifold structure on $\Sigma$,
whose set of orbifold points is contained in the given
marked-point set $\z$.

\vspace{2mm}

Let $(S^2,\z)$ be a genus-zero Riemann surface with k-marked
points $\z=(z_1, \cdots,z_k)$, $p\in X$ any point in an almost
complex orbifold $X$ with $\dim_\C X=n$, and $(V_p,G_p,\pi_p)$ a
local chart at $p$. Then for any k-tuple $\g=(g_1,\cdots,g_k)$
where $g_i\in G_p$, $i=1,\cdots,k$,
there is an orbifold structure on $S^2$ so that
it becomes a complex orbicurve $(S^2,\z,\m)$ where
$\m=(|g_1|,\cdots,|g_k|)$ (here $|g|$ stands for the order of
$g$). If further assuming that $o(\g)=g_1g_2\cdots g_k=1_{G_p}$,
one can construct a rank-n holomorphic orbifold bundle
$E_{p,\g}$ over $(S^2,\z,\m)$, together with an orbifold bundle morphism
$\Phi_{p,\g}:E_{p,\g}\rightarrow TX$ covering the constant map
from $S^2$ to $p\in X$, as we shall see next.

Denote ${\bf 1}_{G_p}=(1_{G_p},\cdots,1_{G_p})$. The case $\g={\bf
1}_{G_p}$ is trivial; we simply take the rank-n trivial
holomorphic bundle over $S^2$. Hence in what follows, we assume that
$\g\neq {\bf 1}_{G_p}$. We recall that the orbifold fundamental group
of $(S^2,\z,\m)$ is given by
$$
\pi^{orb}_1(S^2)=\{\lambda_1,\lambda_2,\cdots,\lambda_k|
   \lambda^{|g_i|}_i=1, \lambda_1\lambda_2\cdots\lambda_k=1\},
$$
where each generator $\lambda_i$ is geometrically represented by a loop
around the marked point $z_i$. We define a homomorphism
$\rho:\pi_1^{orb}(S^2)\rightarrow G_p$ by
sending each $\lambda_i$ to $g_i\in G_p$ (note that we assumed that
$g_1g_2\cdots g_k=1_{G_p}$). There is a closed Riemann surface
$\Sigma$ and a finite group $G$ acting on $\Sigma$ holomorphically,
such that $(\Sigma,G)$ uniformizes $(S^2,\z,\m)$ and
$\pi_1(\Sigma)=\ker\;\rho$ with $G=Im \rho\subset G_p$.
We identify $(TV_p)_p$ with $\C^n$ and denote the rank-n trivial
holomorphic vector bundle on $\Sigma$ by $\underline{\C^n}$.
The representation $G\rightarrow  Aut((TV_p)_p)$ defines a holomorphic action
on the holomorphic vector bundle $\underline{\C}^n$. We take $E_{p,\g}$
to be the corresponding holomorphic orbifold bundle uniformized by
$(\underline{\C}^n,G,\tilde{\pi})$ where
$\tilde{\pi}:\underline{\C}^n\rightarrow \underline{\C}^n/G$ is
the quotient map. There is a natural orbifold bundle morphism $\Phi_{p,\g}:
E_{p,\g}\rightarrow TX$ sending $\Sigma$ to the point $p$.

By the nature of construction, if $\g=(g_1,\cdots,g_k)$ and
$\g^\prime =(g_1^\prime,\cdots,g_k^\prime)$ are conjugate, i.e.,
there is an element $g\in G_p$ such that $g_i^\prime=g^{-1}g_ig$,
then there is an isomorphism $\psi:E_{p,\g}\rightarrow
E_{p,\g^\prime}$ such that $\Phi_{p,\g}=
\Phi_{p,\g^\prime}\circ\psi$.

If there is an isomorphism $\psi:E_{p,\g}\rightarrow
E_{p,\g^\prime}$ such that
$\Phi_{p,\g}=\Phi_{p,\g^\prime}\circ\psi$, then there is a lifting
$\tilde{\psi}:\tilde{E}_{p,\g}\rightarrow \tilde{E}_{p,\g^\prime}$
of $\psi$ and an automorphism $\phi:TV_p\rightarrow TV_p$, such
that $\phi\circ
\tilde{\Phi}_{p,\g}=\tilde{\Phi}_{p,\g^\prime}\circ\tilde{\psi}$.
If $\phi$ is given by the action of an element $g\in G_p$, then we
have $gg_ig^{-1}= g_i^\prime$ for all $i=1,\cdots,k$.
    \vskip 0.1in
    \noindent
    {\bf Lemma 4.2.3: }{\it
Let $E$ be a rank-n holomorphic orbifold bundle over
$(S^2,\z,\m)$ (for some $\m$). Suppose that there is an orbifold
bundle morphism $\Phi:E\rightarrow TX$ covering a constant map
from $S^2$ into $X$. Then there is a $(p,\g)$ such that
$(E,\Phi)=(E_{p,\g},\Phi_{p,\g})$.}

\vspace{3mm}

\noindent{\bf Proof:} Let $E$ be a rank-n holomorphic orbifold
bundle over $(S^2,\z,\m)$ with a morphism $\Phi:E\rightarrow
TX$ covering the constant map to a point $p$ in X. We will find a
$\g$ so that $(E,\Phi)=(E_{p,\g},\Phi_{p,\g})$. For this purpose,
we again consider the uniformizing system $({\Sigma},G,\pi)$ of
$(S^2,\z,\m)$ where ${\Sigma}$ is a closed
Riemann surface with a holomorphic action by a finite group $G$.
Then there is a holomorphic vector bundle $\tilde{E}$ over
${\Sigma}$ with a compatible action of $G$ so that
$(\tilde{E},G)$ uniformizes the holomorphic orbifold bundle $E$.
Moreover, there is a vector bundle morphism
$\tilde{\Phi}:\tilde{E} \rightarrow TV_p$, which is a lifting of
$\Phi$ so that for any $a\in G$, there is a $\tilde{\lambda}(a)$
in $G_p$ such that $\tilde{\Phi}\circ a=
\tilde{\lambda}(a)\circ\Phi$. In fact, $a\rightarrow
\tilde{\lambda}(a)$ defines a homomorphism
$\tilde{\lambda}:G\rightarrow G_p$. Since $\tilde{\Phi}$ covers a
constant map
from ${\Sigma}$ into $V_p$, the holomorphic vector bundle
$\tilde{E}$ is in fact a trivial bundle. Recall that $G$ is the
quotient group of $\pi_1^{orb}(S^2)$ by the normal
subgroup $\pi_1(\Sigma)$. Let
$\lambda$ be the induced homomorphism
$\pi_1^{orb}(S^2)\rightarrow G_p$, and let
$g_i=\lambda(\gamma_i)$. Then we have $g_1g_2\cdots g_k=1_{G_p}$.
We simply define $\g=(g_1,g_2,\cdots,g_k)$. It is easily seen that
$(E,\Phi)=(E_{p,\g},\Phi_{p,\g})$. \hfill $\Box$

\vspace{3mm}

\noindent{\bf Definition 4.2.4: }{\it Given a genus-zero Riemann
surface with k-marked points $(\Sigma,\z)$, where
$\z=(z_1,\cdots,z_k)$, we call each equivalence class
$[E_{p,\g},\Phi_{p,\g}]$ of pair $(E_{p,\g}, \Phi_{p,\g})$ a {\it
ghost map} from $(\Sigma,\z)$ into $X$. A ghost map $[E,\Phi]$
from $(\Sigma,\z)$ is said to be {\it equivalent} to a ghost map
$[E^\prime,\Phi^\prime]$ from $(\Sigma^\prime,\z^\prime)$
($\z^\prime=(z_1^\prime,\cdots,z_k^\prime)$) if there is a
holomorphic orbifold bundle morphism $\tilde{\psi}: E\rightarrow
E^\prime$ covering a biholomorphism $\psi:\Sigma\rightarrow
\Sigma^\prime$ such that $\psi(z_i)=z_i^\prime$ and
$\Phi=\Phi^\prime\circ\tilde{\psi}$. An equivalence class of ghost
maps is called a {\it ghost curve} (with k-marked points). We
denote by $\M_k$ the {\it moduli space of ghost curves with
k-marked points}}. \hfill $\Box$

\vspace{3mm}

As a consequence, we obtain

\vspace{3mm}

\noindent{\bf Proposition 4.2.5: }{\it Let $X$ be an almost
complex orbifold. For any $k\geq 0$, the moduli space of ghost
curves with $k$-marked points $\M_k$ is naturally an almost
complex orbifold. When $k\geq 4$, $\M_k$ can be identified with
$\M_{0,k}\times \widetilde{X}_k^o$, where $\M_{0,k}$ is the
moduli space of genus-zero curve with $k$-marked points. It has a
natural partial compatification $\overline{\M}_k$, which is an
almost complex orbifold and can be identified with
$\overline{\M}_{0,k}\times \widetilde{X}_k^o$, where
$\overline{\M}_{0,k}$ is the Deligne-Mumford compatification of
$\M_{0,k}$.
}

\vspace{3mm}

\noindent{\bf Remarks 4.2.6: }{\it
    \begin{description}
    \item[(i)] The natural partial compatification
$\overline{\M}_k$ of $\M_k$ ($k\geq 4$) can be interpreted
geometrically as adding {\it nodal} ghost curves into $\M_k$.
    \item[(ii)]
  The space $\widetilde{X}_2^o$ is naturally identified with the
graph of the map $I:\widetilde{\Sigma X}\rightarrow
\widetilde{\Sigma X}$ in $\widetilde{\Sigma
X}\times\widetilde{\Sigma X}$, where $I$ is defined by $(p,(g))
\rightarrow (p,(g^{-1}))$.
    \end{description}}

    Next, we   construct a
complex orbifold bundle $E_k$, a kind of obstruction bundle in
nature,  over the moduli space $\M_k$ of ghost curves with
k-marked points. The rank of $E_k$ may vary over different
connected components of $\M_k$. When $k=3$, the restriction of
$E_3$ to each component gives a geometric
    construction of obstruction bundle $E_{(\g)}$ in the last section
    under
    identification $\M_3=\widetilde{X}_3^o$.

Let us consider the space ${\cal{C}}_k$ of all triples
$((\Sigma,\z),E_{p,\g},\Phi_{p,\g})$ where $(\Sigma,\z)$ is a
genus-zero curve with k-marked points $\z=(z_1,\cdots,z_k)$,
$E_{p,\g}$ is a rank-n holomorphic orbifold bundle over $\Sigma$,
and $\Phi_{p,\g}:E_{p,\g}\rightarrow TX$ a morphism covering the
constant map sending $\Sigma$ to the point $p$ in $X$. To each
point $x\in {\cal{C}}_k$ we assign a complex vector space $V_x$,
which is the cokernel of the operator $$ \bar{\partial}:
\Omega^{0,0}(E_{p,\g})\rightarrow \Omega^{0,1}(E_{p,\g}).
\leqno(4.2.11) $$ We introduce an equivalence relation $\sim$
amongst pairs $(x,v)$ where $x\in {\cal{C}}_k$ and $v\in V_x$ as
follows: Let $x=((\Sigma,\z),E_{p,\g},\Phi_{p,\g})$ and
$x^\prime=((\Sigma^\prime,\z^\prime),E_{p^\prime,\g^\prime},\Phi_{p^\prime
,\g^\prime})$, then $(x,v)\sim (x^\prime,v^\prime)$ if there is a
morphism $\tilde{\psi}:E_{p,\g}\rightarrow E_{p^\prime,\g^\prime}$
such that $\Phi_{p,\g}=\Phi_{p^\prime,\g^\prime}\circ\tilde{\psi}$
and $\tilde{\psi}$ covers a biholomorphism
$\psi:\Sigma\rightarrow\Sigma^\prime$ satisfying
$\psi(\z)=\z^\prime$ (as ordered sets),  and
$v^\prime=\psi_\ast(v)$ where $\psi_\ast:V_x\rightarrow
V_{x^\prime}$ is induced by $\tilde{\psi}$. We define $E_k$ to be
the quotient space of all $(x,v)$ under this equivalence relation.
There is obviously a surjective map $pr:E_k\rightarrow \M_k$
induced by the projection $(x,v)\rightarrow x$.

\vspace{3mm}

\noindent{\bf Lemma 4.2.7:} The space $E_k$ can be given a
topology such that $pr:E_k\rightarrow \M_k$ is a complex orbifold
bundle over $\M_k$.

\vspace{3mm}

\noindent{\bf Proof:} First we show that the dimension of $V_x$ is
a local constant function of the equivalence class $[x]$ in
$\M_k$. Recall a neighborhood of $[x]$ in $\M_k$ is given by
${\cal{O}}\times V_p^\g/C(\g)$ where ${\cal{O}}$ is a neighborhood
of the genus-zero curve with k-marked points $(\Sigma,\z)$ in the
moduli space $\M_{0,k}$. In fact, we will show that the kernel of
$(4.2.11)$ is identified with $(TV_p^\g)_p$, whose dimension is a
local constant. Then it follows that $\dim V_x$ is locally
constant as the dimension of cokernel of $(4.2.11)$, since by
Proposition 4.2.2, the index of $(4.2.11)$ is locally constant.

For the identification of the kernel of $(4.2.11)$, recall that
the holomorphic orbifold bundle $E_{p,\g}$ over the genus-zero
curve $\Sigma$ is uniformized by the trivial holomorphic vector
bundle $\underline{\C}^n$ over a Riemann surface $\tilde{\Sigma}$
with a holomorphic action of a finite group $G$. Hence the kernel
of $(4.2.11)$ is identified with the $G$-invariant holomorphic
sections of the trivial bundle $\underline{\C}^n$, which are
constant sections invariant under $G$. Through morphism
$\Phi_{p,\g}: E_{p,\g}\rightarrow TX$, the kernel of
$\bar{\partial}$ is then identified with $(TV_p^\g)_p$.

Recall that the moduli space $\M_{0,k}$ is a smooth complex
manifold. Let ${\cal{O}}$ be a neighborhood of $(\Sigma_0,\z_0)$
in $\M_{0,k}$. Then a neighborhood of
$[x_0]=[(\Sigma_0,\z_0),E_{p,\g},\Phi_{p,\g}]$ in $\M_k$ is
uniformized by $({\cal{O}}\times V_p^\g, C(\g))$ (cf. Lemma 4.1.1).
More precisely, to any $((\Sigma,\z),y)\in {\cal{O}}\times
V_p^\g$, we associate a rank-n holomorphic orbifold bundle over
$(\Sigma,\z)$ as follows: Let $q=\pi_p(y)\in U_p$, then the pair
$(y,\g)$ canonically determines a $\h_y\in G_q\times\cdots\times
G_q$, and there is a canonically constructed holomorphic orbifold
bundle $E_{q,\h_y}$ over $(\Sigma,\z)$ with morphism
$\Phi_{q,\h_y}:E_{q,\h_y}\rightarrow TX$ covering the constant map
to $q$. Hence we have a family of holomorphic orbibundles over
genus-zero curve with k-marked points, which are parametrized by
${\cal{O}}\times V_p^\g$. Moreover, it depends on the parameter in
${\cal{O}}$ holomorphically and the action of $C(\g)$ on $V_p^\g$
coincides with the equivalence relation between the pairs of
holomorphic orbifold bundle and morphism
$(E_{q,\h_y},\Phi_{q,\h_y})$. Now we put a Kahler metric on each
genus-zero curve in ${\cal{O}}$ which is compatible to the complex
structure and depends smoothly on the parameter in ${\cal{O}}$,
and we also put a hermitian metric on $X$. Then we have a family
of first order elliptic operators depending smoothly on the
parameters in ${\cal{O}}\times V_p^\g$: $$
\bar{\partial}^\ast:\Omega^{0,1}(E_{q,\h_y})\rightarrow
\Omega^{0,0}(E_{q,\h_y}) $$ and whose kernel gives rise to a
complex vector bundle $E_{x_0}$ over ${\cal{O}}\times V_p^\g$. The
finite group $C(\g)$ naturally acts on the complex vector bundle
which coincides with the equivalence relation amongst the pairs
$(x,v)$ where $x\in{\cal{C}}_k$ and $v\in V_x$. Hence
$(E_{x_0},C(\g))$ is a uniformizing system for
$pr^{-1}({\cal{O}}\times V_p^\g/C(\g))$, which fits together to
give an orbifold bundle structure for $pr:E_k\rightarrow \M_k$.
\hfill $\Box$

\vspace{3mm}

\noindent{\bf Remark 4.2.8:} Recall that each holomorphic orbifold
bundle $E_{p,\g}$ over $(S^2,\z,\m)$ can be uniformized by a trivial
holomorphic vector bundle $\underline{\C}^n$ over a Riemann
surface ${\Sigma}$ with a holomorphic group action by $G$.
Hence each element $\xi$ in the kernel of $$
\bar{\partial}^\ast:\Omega^{0,1}(E_{p,\g})\rightarrow
\Omega^{0,0}(E_{p,\g}) $$ can be identified with a $G$-invariant
harmonic $(0,1)$-form on ${\Sigma}$ with value in $(TV_p)_p$
(here we identify each fiber of $\underline{\C}^n$ with $(TV_p)_p$
through the morphism $\Phi_{p,\g}$), i.e., $\xi=w\otimes\alpha$
where $w\in (TV_p)_p$, $\alpha$ is a harmonic $(0,1)$-form on
$\tilde{\Sigma}$, and $\xi$ is $G$-invariant. Therefore, when
$k=3$, it agrees with $E_{(\g)}$. We observe that with respect to the
taken hermitian metric on $X$, $w\in (TV_p)_p$ must lie in the
orthogonal complement of $(TV_p^\g)_p$ in $(TV_p)_p$. This is
because: For any $u\in (TV_p^\g)_p$ and a harmonic $(0,1)$-form
$\beta$ on ${\Sigma}$, if $u\otimes\beta$ is $G$-invariant,
then $\beta$ is $G$-invariant too, which means that $\beta$
descents to a harmonic $(0,1)$-form on $S^2$, and
$\beta$ must be identically zero.
\hfill $\Box$

\vspace{3mm}

    Recall the cup product is defined by equation

$$ <\alpha_1\cup_{orb}\alpha_2, \gamma>_{orb}=
 (\int_{\M_3}^{orb} e_1^\ast(\alpha_1)\cup
e_2^\ast(\alpha_2)\cup e_3^\ast(\gamma)\cup e(E_3)), $$ where
$e(E_3)$ is the Euler form of the complex orbifold bundle $E_3$
over $\M_3$ and $\gamma\in H^{*}_{orb,c}(X)$.

\vspace{3mm}

    We take a basis $\{e_j\}, \{e^o_k\}$ of the total orbifold
    cohomology group $H^*_{orb}(X), H^*_{orb,c}(X)$ such that each
    $e_j, e^o_k$ is of
    homogeneous degree. Let $<e_j, e^o_k>_{orb}=a_{jk}$ be the Poincare pairing
    matrix and $(a^{jk})$ be the inverse. It is easy to check that the Poincare
    dual of graph of $I$ in $\widetilde{\Sigma}^2$ can be written as
    $\sum_{j,k}a^{jk}e_j\otimes e^o_k$. Then,
$$ \alpha_1\cup_{orb}\alpha_2= \sum_{j,k}e_j
a^{kj}(\int_{\M_3}^{orb} e_1^\ast(\alpha_1)\cup e_2
^\ast(\alpha_2)\cup e_3^\ast(e_k^o)\cup e(E_3)). \leqno(4.2.12) $$

\noindent{\bf Proof of Theorem 4.1.5:} We postpone the proof of
associativity of $\cup_{orb}$ to the next subsection.

We first show that if $\alpha_1\in H^p_{orb}(X)$ and
$\alpha_2\in H^q_{orb}(X)$, then $\alpha_1\cup_{orb}\alpha_2$
is in $H^{p+q}_{orb}(X)$. For the integral in (4.2.12) to be
nonzero, $$
\deg(e_1^\ast(\alpha_1))+\deg(e_2^\ast(\alpha_2))+\deg(e_3^\ast(e_k^o))+
\deg(e(E_3))=2\dim_\C\M_3.\leqno(4.2.13) $$ Here $\deg$ stands for
the degree of a cohomology class without degree shifting. The
degree of Euler class $e(E_3)$ is equal to the dimension of
cokernel of $(4.2.11)$, which by index formula (cf. Proposition
4.2.2) equals $2\dim_\C \M_3^{(i)}-(2n-2\sum_{j=1}^3
\iota(p,g_j))$ on a connected component $\M_3^{(i)}$ containing
point $(p,(\g))$ where $\g=(g_1,g_2,g_3)$. Hence $(4.2.13)$ becomes
$$ \deg(\alpha_1)+deg(\alpha_2)+deg(e_k^o)+2\sum_{j=1}^3
\iota(p,g_j)=2n,\leqno(4.2.14) $$ from which it is easily seen that
$\alpha_1\cup_{orb}\alpha_2$ is in $H^{p+q}_{orb}(X)$.

Next we show that $e_X^0$ is a unit with respect to $\cup_{orb}$,
i.e., $\alpha\cup_{orb} e^0_X=e^0_X\cup_{orb}\alpha=\alpha$. First
observe that there are connected components of $\M_3$ consisting
of points $(p,(\g))$ for which $\g=(g_1,g_2,g_3)$ satisfies the
condition that one of the $g_i$ is $1_{G_p}$. Over these
components the Euler class $e(E_3)=1$ in the $0^{th}$ cohomology
group since $(4.2.11)$ has zero cokernel. Let $\alpha\in
H^\ast(X_{(g)})$. Then
$e_1^\ast(\alpha)\cup e_2^\ast(e_X^0)\cup e_3^\ast(e_k^o)$ is
non-zero only on the connected component of $\M_3$  which is the
image of the embedding $X_{(g)}\rightarrow\M_3$ given by
$(p,(g)_{G_p})\rightarrow (p,((g,1_{G_p},g^{-1})))$ and $e_k^o$ must
be in $H_c^\ast(X_{(g^{-1})})$. Moreover, we have

\begin{eqnarray*}
\alpha\cup_{orb} e^0_X & := &
\sum_{j,k}(\int_{\M_3}^{orb} e_1^\ast(\alpha)\cup e_2^\ast(e_X^0)\cup
e_3^\ast(e_k^o)\cup e(E_3))a^{kj}e_j\\
                     & =  & \sum_{j,k}(\int_{X_{(g)}}^{orb} \alpha\cup
                            I^\ast(e_k^o))a^{kj}e_j\\
                     & =  & \alpha
\end{eqnarray*}
Similarly, we can prove that $e^0_X\cup_{orb}\alpha=\alpha$.

Now we consider the case $\cup_{orb}: H^d_{orb}(X)\times
H^{2n-d}_{orb}(X) \rightarrow H^{2n}_{orb}(X)=H^{2n}(X)$.
Let $\alpha\in H^d_{orb}(X)$ and $\beta\in
H^{2n-d}_{orb}(X)$, then
$e_1^\ast(\alpha)\cup e_2^\ast(\beta)\cup e_3^\ast (e_X^0)$ is
non-zero only on those connected components of $\M_3$ which are
images under embedding $\widetilde{X}\rightarrow \M_3$ given by
$(p,(g))\rightarrow (p,((g,g^{-1},1_{G_p})))$, and if $\alpha$ is
in $H^\ast(X_{(g)})$, $\beta$ must be in
$H^\ast(X_{(g^{-1})})$. Moreover, let $e_X^{2n}$ be the
generator in $H^{2n}(X)$ such that $e_X^{2n}\cdot [X]=1$, then
we have

\begin{eqnarray*}
\alpha\cup_{orb}\beta & := &
\sum_{j,k}(\int_{\M_3}^{orb} e_1^\ast(\alpha)\cup e_2^\ast(\beta)
\cup e_3^\ast (e_k^o)\cup e(E_3))a^{kj}e_j\\
                   & =  & (\int_{\M_3}^{orb} e_1^\ast(\alpha)\cup e_2^\ast
(\beta)\cup e_3^\ast (e_X^0)\cup e(E_3))\cdot e_X^{2n}\\
                     & =  & (\int_{\widetilde{X}}^{orb} \alpha\cup
                             I^\ast(\beta))\cdot e_X^{2n}\\
                     & = &<\alpha, \beta>_{orb} e_X^{2n}
\end{eqnarray*}
from which we see that $\int_X\alpha\cup_{orb}\beta
=<\alpha,\beta>_{orb}$.

The rest of the assertions are obvious. \hfill $\Box$

\subsection{Proof of associativity}

In this subsection, we give a proof of associativity of the
orbifold cup products $\cup_{orb}$ defined in the last subsection.
We will only present the proof for the orbifold cohomology groups
$H^\ast_{orb}(X)$. The proof for orbifold Dolbeault cohomology
is the same. We leave it to readers.

\vspace{3mm}

Recall the moduli space of ghost curves with k-marked points
$\M_k$ for $k\geq 4$ can be identified with $\M_{0,k}\times
\widetilde{X}_k^o$ which admits a natural partial
compatification $\overline{\M}_{0,k}\times \widetilde{X}_k^o$ by
adding nodal ghost curves. We will first give a detailed
analysis on this for the case when $k=4$.

Let $\Delta$ be the graph of map $I:\widetilde{\Sigma
X}\rightarrow \widetilde{\Sigma X}$ in $\widetilde{\Sigma X}\times
\widetilde{\Sigma X}$ given by $I: (p,(g))\rightarrow
(p,(g^{-1}))$. To obtain the orbifold structure, one can view
$\Delta$ as orbifold fiber product of identify map and $I$, which
has an induced orbifold structure since both identify and $I$ are
so called "good map" (see \cite{CR}).  Consider map $\Lambda:
\tilde{X}^o_3\times\tilde{X}^o_3\rightarrow \widetilde{\Sigma X}
\times\widetilde{\Sigma X}$ given by $((p,(\g)),
(q,(\h)))\rightarrow ((p,(g_3)),(q,(h_1)))$. We wish to consider
the preimage of $\Delta$.
    \vskip 0.1in
    \noindent
    {\bf Remark: }{\it Suppose that we have two maps
    $$f: X\rightarrow Z, g: Y\rightarrow Z.$$
    In general, ordinary fiber product $X\times_Z Y$ may not have
    a natural orbifold structure. The correct formulation is to use
    "good map" introduced in \cite{CR}. If $f,g$ are good maps, there is
    a canonical orbifold fiber product (still denoted by $X\times_Z Y$)
    obtained by taking fiber product on uniformizing system. It has
    an induced orbifold structure and there are good map projection to both $X, Y$
    to make appropriate diagram to commute. However, as a set, such an
    orbifold fiber product is not usual fiber product. Throughout this paper,
    we will use $X\times_Z Y$ to denote orbifold fiber product only.}
    \vskip 0.1in
    It is clear that the pre-image of $\Delta$ can be viewed as fiber product of
    $$e_3, I\circ e_1: \tilde{X}^0_3\rightarrow \widetilde{X}.$$
    Then, we define the pre-image $\Lambda^{-1}(\Delta)$ as orbifold fiber product
    of $e_3, , I\circ e_1$. It is easy to check that $\Lambda^{-1}(\Delta)=\tilde{X}^o_4$.
     Next, we describe
explicitly the compatification $\overline{\M}_4$ of $\M_4$.

Recall the moduli space of genus-zero curves with 4-marked points
$\M_{0,4}$ can be identified with ${\bf P}^1\setminus
\{0,1,\infty\}$ by fixing the first three marked points to be
$\{0,1,\infty\}$. The Deligne-Mumford compactification
$\overline{\M}_{0,4}$ is then identified with ${\bf P}^1$ where
each point of $\{0,1,\infty\}$ corresponds to a nodal curve
obtained as the last marked point is running into this point. It
is easy to see that part of the compatification $\overline{\M}_4$
by adding a copy of $\widetilde{X}_4^o$ at $\infty \in
\overline{\M}_{0,4}={\bf P}^1$ where intuitively we associate
$(g_1g_2)^{-1}, g_1g_2$ at nodal point.  In the same way, the
compatification at $0$ is by adding a copy of $\widetilde{X}_4^o$
where we associate $(g_1g_4)^{-1}, g_1g_4$ at nodal point,   and
at $1$ by associating $(g_1g_3)^{-1}, g_1g_3$ at nodal point.

    Next, we define an orbifold bundle to measure the failure of
    transversality of $\Lambda$ to $\Delta$.
\vspace{3mm}

\noindent{\bf Definition 4.3.1:} {\it We define a complex orbifold
bundle $\nu$ over $\Lambda^{-1}(\Delta)_{(g_1,g_2,g_3,g_4)}$ as
follows: over each uniformizing system $(V^\g_p,C(\g))$ of
$\Lambda^{-1}(\Delta_{(\g)})$, where $\g=(g_1,g_2,g_3,g_4)$, we
regard $V^\g_p$ as the intersection of $V_p^{g_1}\cap V_p^{g_2}$
with $V_p^{g_3}\cap V_p^{g_4}$ in $V_p^g$ where $g=(g_1g_2)^{-1}$.
We define $\nu$ to be the complex orbifold bundle over
$\Lambda^{-1}(\Delta)$ whose fiber is the orthogonal complement of
$V_p^{g_1}\cap V_p^{g_2} + V_p^{g_3}\cap V_p^{g_4}$ in $V_p^g$. }

\vspace{3mm}

The associativity is based on the following

\vspace{3mm}

\noindent{\bf Lemma 4.3.2: }{\it The complex orbifold bundle
$pr:E_4\rightarrow \M_4$ can be extended over the compatification
$\overline{\M}_4$, denoted by
$\bar{pr}:\overline{E}_4\rightarrow\overline{\M}_4$, such that
$\overline{E}_4|_{\{\ast\}\times \widetilde{X}_4^o}= (E_3\oplus
E_3)|_{\Lambda^{-1}(\Delta)}\oplus\nu$ under the above identification, where
$\{\ast\}$ represents a point in
$\{0,1,\infty\}\subset\overline{\M}_{0,4}$.}

\vspace{3mm}

\noindent{\bf Proof:} We fix an identification of infinite
cylinder $\R\times S^1$ with $\C^\ast \setminus\{0\}$ via the
biholomorphism defined by $t+is\rightarrow e^{-(t+is)}$ where
$t\in\R$ and $s\in S^1=\R/2\pi\Z$. Through this identification, we
regard a punctured Riemann surface as a Riemann surface with
cylindrical ends. A neighborhood of a point $\ast\in
\{0,1,\infty\}\subset\overline{\M}_{0,4}$,  as a family of
isomorphism classes of genus-zero curves with 4-marked points, can
be described by a family of curves $(\Sigma_{r,\theta},\z)$
obtained by gluing of two genus-zero curves with a cylindrical end
and two marked points on each, parametrized by $(r,\theta)$ where
$0\leq r\leq r_0$ and $\theta\in S^1$, as we glue the two curves
by self-biholomorphisms of $(-\ln r,-3\ln r)\times S^1$ defined by
$(t,s)\rightarrow (-4\ln r-t,-(s+\theta))$ ($r=0$ represents the
nodal curve $\ast$). Likewise, thinking of points in $\M_4$ as
equivalence classes of triples
$((\Sigma,\z),E_{p,\g},\Phi_{p,\g})$ where $(\Sigma,\z)$ is a
genus-zero curve of 4-marked points $\z$, a neighborhood of
$\{\ast\}\times (X\sqcup\widetilde{X}_4^o)$ in $\overline{\M}_4$
are described by a family of holomorphic orbifold bundles on
$(\Sigma_{r,\theta},\z)$ with morphisms obtained by gluing two
holomorphic orbifold bundles on genus-zero curves with two marked
points and one cylindrical end on each. We denote them by
$(E_{r,\theta},\Phi_{r,\theta})$.

The key is to construct a family of isomorphisms of complex
orbifold bundle $$ \Psi_{r,\theta}: E_3\oplus
E_3\oplus\nu|_{\Lambda^{-1}(\Delta)}\rightarrow E_4 $$ for
$(r,\theta)\in (0,r_0)\times S^1$. Recall the fiber of $E_3$ and
$E_4$ is given by kernels of the $\bar{\partial}^\ast$ operators.
In fact, $\Psi_{r,\theta}$ are given by gluing maps of kernels of
$\bar{\partial}^\ast$ operators.

More precisely, suppose
$((\Sigma_{r,\theta},z),E_{r,\theta},\Phi_{r,\theta})$ are
obtained by gluing $((\Sigma_1,\z_1), E_{p,\g},\Phi_{p,\g})$ and
$((\Sigma_2,\z_2), E_{p,\h},\Phi_{p,\h})$ where $\g=(g_1,g_2,g)$
and $\h=(g^{-1},h_2,h_3)$. Let $m=|g|$. Then $E_{r,\theta}|_{(-\ln
r,-3\ln r)\times S^1}$ is uniformized by $(-\frac{\ln
r}{m},-\frac{3\ln r}{m})\times S^1\times TV_p$ with an obvious
action by $\Z_m=\langle g\rangle$.

Let $\xi_1\in \Omega^{0,1}(E_{p,\g})$, $\xi_2\in
\Omega^{0,1}(E_{p,\h})$ such that $\bar{\partial}^\ast\xi_i=0$ for
$i=1,2$. On the cylindrical end, if we fix the local coframe
$d(t+is)$, then each $\xi_i$ is a $TV_p$-valued, exponentially
decaying holomorphic function on the cylindrical end. We fix a
cut-off function $\rho(t)$ such that $\rho(t)\equiv 1$ for $t\leq
0$ and $\rho(t)\equiv 0$ for $t\geq 1$. We define the gluing of
$\xi_1$ and $\xi_2$, which is a section of
$\Omega^{0,1}(E_{r,\theta})$ and denoted by $\xi_1 \#\xi_2$, by $$
\xi_1 \#\xi_2 = \rho(-2\ln r+t)\xi_1+(1-\rho(-2\ln r+t))\xi_2 $$
on the cylindrical end. Let $\Psi_{r,\theta}(\xi_1,\xi_2)$ be the
$L^2$-projection of $\xi_1\#\xi_2$ onto $\ker\bar{\partial}^\ast$,
then the difference $\eta=\xi_1 \#\xi_2
-\Psi_{r,\theta}(\xi_1,\xi_2)$ satisfies the estimate
$||\bar{\partial}^\ast\eta||_{L^2} \leq
Cr^{\delta}(||\xi_1||+||\xi_2||)$ for some
$\delta=\delta(\xi_1,\xi_2)>0$. Hence $||\eta||_{L^2}\leq C|\ln
r|r^\delta(||\xi_1||+||\xi_2||)$ (cf. \cite{Ch}), from which it
follows that for small enough $r$, $\Psi_{r,\theta}$ is an
injective linear map.

Now given any $\xi\in V_p^g$ which is orthogonal to both
$V_p^{g_1}\cap V_p^{g_2}$ and $V_p^{g_3}\cap V_p^{g_4}$, we define
$\Psi_{r,\theta}(\xi)$ as follows: fixing a cut-off function, we
construct a section $u_\xi$ over the cylindrical neck $(-\ln
r,-3\ln r)\times S^1$ with support in $(-\ln r+1,-3\ln r-1)\times
S^1$ and equals $\xi$ on $(-\ln r+2,-3\ln r-2)\times S^1$. We
write $\bar{\partial}^\ast u_\xi=v_{\xi,1}+v_{\xi,2}$ where
$v_{\xi,1}$ is supported in $(-\ln r+1,-\ln r+2)\times S^1$ and
$v_{\xi,2}$ in $(-3\ln r-2,-3\ln r-1)\times S^1$. Since $\xi$ is
orthogonal to both $V_p^{g_1}\cap V_p^{g_2}$ and $V_p^{g_3}\cap
V_p^{g_4}$, we can arrange so that $v_{\xi,1}$ is $L^2$-orthogonal
to $V_p^{g_1}\cap V_p^{g_2}\cap V_p^g$ and $v_{\xi,2}$ is
$L^2$-orthogonal to $V_p^{g^{-1}}\cap V_p^{g_3}\cap V_p^{g_4}$,
which are the kernels of the $\bar{\partial}$ operators on
$\Sigma_1$ and $\Sigma_2$ acting on sections of $E_{p,\g}$ and
$E_{p,\h}$ respectively. Hence there exist $\alpha_1\in
\Omega^{0,1}(E_{p,\g})$ and $\alpha_2\in \Omega^{0,1}(E_{p,\h})$
such that $\bar{\partial}^\ast \alpha_i =v_{\xi,i}$ and $\alpha_i$
are $L^2$-orthogonal to the kernels of the $\bar{\partial}^\ast$
operators respectively. We define $\Psi_{r,\theta}(\xi)$ to be the
$L^2$-orthogonal projection of $u_\xi-\alpha_1\#\alpha_2$ onto
$\ker\,\bar{\partial}^\ast$, then $\Psi_{r,\theta}(\xi)$ is linear
on $\xi$. On the other hand, observe that
$||\bar{\partial}^\ast(u_\xi-\alpha_1\#\alpha_2)||_{L^2}\leq
Cr^\delta ||\xi||$ for some $\delta>0$, if let $\eta$ be the
difference of $\Psi_{r,\theta}(\xi)$ and
$u_\xi-\alpha_1\#\alpha_2$, then $||\eta||_{L^2}\leq C|\ln
r|r^\delta||\xi||$ (cf. \cite{Ch}), from which we see that for
sufficiently small $r>0$, $\Psi_{r,\theta}(\xi)\neq 0$ if $\xi\neq
0$.

Hence we construct a family of injective morphisms $$
\Psi_{r,\theta}: E_3\oplus
E_3\oplus\nu|_{\Lambda^{-1}(\Delta)}\rightarrow E_4 $$ for
$(r,\theta)\in (0,r_0)\times S^1$. We will show next that each
$\Psi_ {r,\theta}$ is actually an isomorphism.

We denote by $\bar{\partial}_i$ the $\bar{\partial}$ operator on
$\Sigma_i$, and $\bar{\partial}_{r,\theta}$ the $\bar{\partial}$
operator on $\Sigma_{r,\theta}$. Then index formula tells us that
(cf. Proposition 4.2.2)
\begin{eqnarray*}
index\,\bar{\partial}_1 & = & n-\sum_{j=1}^3\iota(p,g_j),\\
index\,\bar{\partial}_2 & = & n-\sum_{j=1}^3\iota(p,h_j),\\
index\,\bar{\partial}_{r,\theta} & = &
n-(\iota(p,g_1)+\iota(p,g_2)+ \iota(p,h_2)+\iota(p,h_3)),
\end{eqnarray*}
from which we see that
$index\,\bar{\partial}_1+index\,\bar{\partial}_2
=index\,\bar{\partial}_{r,\theta}+\dim_\C V_p^g$. Since
$\dim\ker\bar{\partial}_1+\dim\ker\bar{\partial}_2=
\dim\ker\bar{\partial}_{r,\theta}+\dim_\C V_p^g-rank\, \nu$, we
have $$ \dim coker\bar{\partial}_1 + \dim coker\bar{\partial}_2 +
rank\, \nu=\dim coker\bar{\partial}_{r,\theta}. $$ Hence
$\Psi_{r,\theta}$ is an isomorphism for each $(r,\theta)$. \hfill
$\Box$

\vspace{3mm}
    Before we prove the associativity, let's review some of basic
    construction of smooth manifold and its orbifold analogue.
    Recall that if $Z\subset X$ is a submanifold, then Poincare
    dual of $Z$ can be constructed by Thom form of normal bundle
    $N_Z$ via the natural identification between normal bundle and
    tubuler neighborhood of $Z$. Here, Thom form $\Theta_Z$ is a close form
    such that its restriction on each fiber is a compact supported
    form of top degree with volume one. In orbifold category, the
    same is true provided that we interpret ``suborbifold''
    correctly. Here, a suborbifold is a good map $f: Z\rightarrow
    X$ such that locally, $f$ can be lifted to a $G$-invariant
    embedding  to ``general'' uniformizing system
    $\tilde{f}: (U_Z, G, \pi_Z)\rightarrow
    (U_X, G, \pi_X)$. Here, ``general'' means that $U_Z, U_X$ could
    be disconnected. For example, orbifold fiber product
    $\Lambda^{-1}(\Delta)$ is a suborbifold of
    $\tilde{X}^o_3\times \tilde{X}^0_3$. It is clear that Poincare
    dual of $Z$ can be represented by Thom class of normal bundle
    $Z$.

    \vskip 0.1in
    \noindent
    {\bf Proposition 4.3.4: }{\it Choose a basis $\{e_j\}, \{e^o_k\}$ of the total orbifold
    cohomology group $H^*_{orb}(X), H^*_{orb,c}(X)$ such that each
    $e_j, e^o_k$ is of
    homogeneous degree. Let $<e_j, e^o_k>_{orb}=a_{jk}$ be the Poincare pairing
    matrix and $(a^{jk})$ be the inverse. Then,
    $$\int_{(\widetilde{X}_4^o)_{(\g)}}^{orb}
e_1^\ast(\alpha_1)\cup e_2^\ast(\alpha_2) \cup e_3^\ast(\alpha_3)
\cup e_4^\ast(e_l^o)\cup e(E_4) $$
    $$=\sum_{j,k}
(\int_{\tilde{X}^o_3}^{orb} e_1^\ast(\alpha_1)\cup e_2^\ast(\alpha_2)
\cup e_3^\ast(e_k^o)\cup e(E_3)) \cdot (\int_{\tilde{X}^o_3}^{orb}
e_1^\ast(e_j)\cup e_2^\ast(\alpha_3) \cup e_3^\ast(e_l^o)\cup
e(E_3))\cdot a^{kj} $$
    }
    \vskip 0.1in
    {\bf Proof: } Key observation is
    $\Lambda^*N_{\Delta}=N_{\Lambda^{-1}(\Delta)}\oplus \nu.$
    Hence,
    $\Lambda^*\Theta_{\Delta}=\Theta_{\Lambda^{-1}(\Delta)}
    \cup \Theta_{\nu}$.
    $$\begin{array}{ll}
    &\int_{\widetilde{X}_4^o}^{orb}
e_1^\ast(\alpha_1)\cup e_2^\ast(\alpha_2) \cup e_3^\ast(\alpha_3)
\cup e_4^\ast(e_l^o)\cup e(E_4)\\
=&\int_{\Lambda^{-1}(\Delta)}^{orb} e_1^\ast(\alpha_1)\cup
    e_2^\ast(\alpha_2) \cup e_3^\ast(\alpha_3) \cup
    e_4^\ast(e_l^o)\cup e(E_3)\cup e(E_3)\cup e(\nu)\\
    =&\int_{\tilde{X}^o_3\times \tilde{X}^o_3}^{orb}
e_1^\ast(\alpha_1)\cup e_2^\ast(\alpha_2)
    \cup e_3^\ast(\alpha_3) \cup e_4^\ast(e_l^o)\cup e(E_3)\cup
    e(E_3)\cup \Theta_{\nu}\cup \Theta_{\Lambda^{-1}(\Delta)}\\
    =&\int_{\tilde{X}^o_3\times \tilde{X}^o_3}^{orb}
e_1^\ast(\alpha_1)\cup e_2^\ast(\alpha_2)
    \cup e_3^\ast(\alpha_3) \cup e_4^\ast(e_l^o)\cup e(E_3)\cup e(E_3)
    \cup \Lambda^\ast \Theta_{\Delta}\\
 =&\sum_{j,k}
(\int_{\tilde{X}^o_3}^{orb} e_1^\ast(\alpha_1)\cup e_2^\ast(\alpha_2)
\cup e_3^\ast(e_k^o)\cup e(E_3))) \cdot (\int_{\tilde{X}^o_3}^{orb}
e_1^\ast(e_j)\cup e_2^\ast(\alpha_3) \cup e_3^\ast(e_l^o)\cup
e(E_3))\cdot a^{kj}
    \end{array}$$

Now we are ready to prove

\vspace{3mm}

\noindent{\bf Proposition 4.3.5: }{\it The cup product
$\cup_{orb}$ is associative, i.e., for any $\alpha_i$, $i=1,2,3$,
we have $$ (\alpha_1\cup_{orb}\alpha_2)\cup_{orb}
\alpha_3=\alpha_1\cup_{orb}(\alpha_2\cup_{orb}\alpha_3). $$}

\vspace{3mm}

\noindent{\bf Proof:} By definition of cup product $\cup_{orb}$,
we have $(\alpha_1\cup_{orb}\alpha_2)\cup_{orb}\alpha_3$ equals $$
\sum_{j,k,l,s} (\int_{\M_3}^{orb} e_1^\ast(\alpha_1)\cup
e_2^\ast(\alpha_2) \cup e_3^\ast(e_k^o)\cup e(E_3)) \cdot
(\int_{\M_3}^{orb} e_1^\ast(e_j)\cup e_2^\ast(\alpha_3) \cup
e_3^\ast(e_l^o)\cup e(E_3)) \cdot a^{kj} a^{ls}e_s $$ and
$\alpha_1\cup_{orb}(\alpha_2\cup_{orb}\alpha_3)$ equals $$
    \sum_{j,k,l,s} (\int_{\M_3}^{orb} e_1^\ast(\alpha_1)\cup e_2^\ast(e_j)
\cup e_3^\ast(e_l^o)\cup e(E_3)) \cdot (\int_{\M_3}^{orb}
e_1^\ast(\alpha_2)\cup e_2^\ast(\alpha_3) \cup e_3^\ast(e_k^o)\cup
e(E_3)) \cdot a^{kj} a^{ls}e_s. $$ By Proposition 4.3.4, $$
\sum_{j,k} (\int_{\M_3}^{orb} e_1^\ast(\alpha_1)\cup e_2^\ast(\alpha_2)
\cup e_3^\ast(e_k^o)\cup e(E_3)) \cdot (\int_{\M_3}^{orb}
e_1^\ast(e_j)\cup e_2^\ast(\alpha_3) \cup e_3^\ast(e_l^o)\cup
e(E_3))\cdot a^{kj} $$ equals $$
\int_{\{\infty\}\times\widetilde{X}_4^o}^{orb} e_1^\ast(\alpha_1)\cup
e_2^\ast(\alpha_2) \cup e_3^\ast(\alpha_3) \cup
e_4^\ast(e_l^o)\cup e(E_4), $$ and $$ \sum_{j,k} (\int_{\M_3}^{orb}
e_1^\ast(\alpha_1)\cup e_2^\ast(e_j) \cup e_3^\ast(e_l^o)\cup
e(E_3)) \cdot (\int_{\M_3}^{orb} e_1^\ast(\alpha_2)\cup
e_2^\ast(\alpha_3) \cup e_3^\ast(e_k^o)\cup e(E_3)) \cdot a^{kj}
$$ equals $$ \int_{\{0\}\times\widetilde{X}_4^o}^{orb}
e_1^\ast(\alpha_1)\cup e_2^\ast(\alpha_2) \cup
e_3^\ast(\alpha_3)\cup e_4^\ast(e_l^o)\cup e(E_4). $$ Hence
$(\alpha_1\cup_{orb}\alpha_2)\cup_{orb}\alpha_3
=\alpha_1\cup_{orb} (\alpha_2\cup_{orb}\alpha_3)$. \hfill $\Box$

\vspace{3mm}

\section{Examples}
     In general, it is easy to compute orbifold cohomology once we
     know the action of local group.
     \vskip 0.1in
     \noindent
     {\bf Example 5.1-Kummer surface: } Consider Kummer surface
     $X=T^4/\tau$, where $\tau$ is the involution $x\rightarrow -x$.
     $\tau$ has 16 fixed points, which give 16 twisted
     sectors. It is easily seen that $\iota_{(\tau)}=1$. Hence, we should shift
     the cohomology classes of a twisted sector by $2$ to obtain
     16 degree two classes in  orbifold cohomology. The cohomology
     classes of nontwisted sector come from invariant cohomology
     classes of $T^4$. It is easy to compute that $H^0(X, \R), H^4(X,
     \R)$ has dimension one and $H^2(X, \R)$ has dimension 6.
     Hence, we obtain
     $$b_0^{orb}=b^{orb}_4=1, b^{orb}_1=b^{orb}_3=0,
     b^{orb}_2=22.$$
     Note that orbifold cohomology group of $T^4/\tau$ is
     isomorphic to ordinary cohomology of $K3$-surface, which is
     the the crepant resolution of $T^4/\tau$. However, it is easy
     to compute that Poincare pairing of $H^*_{orb}(T^4/\tau, \R)$
     is different from Poincare pairing of $K3$-surface. We leave
     it to readers

  \vskip 0.1in

     \noindent{\bf Example 5.2-Borcea-Voisin threefold: } An important class of
     Calabi-Yau 3-folds due to Borcea-Voisin is constructed as
     follows: Let $E$ be an elliptic curve with an involution $\tau$
     and $S$ be a $K3$-surface with an involution $\sigma$ acting by
     $(-1)$ on $H^{2,0}(S)$. Then,
     $\tau\times \sigma$ is an involution of $E\times S$, and
     $X=E\times S/<\tau\times \sigma>$ is a Calabi-Yau orbifold. The crepant
     resolution $\widetilde{X}$ of $X$ is a smooth Calabi-Yau 3-fold.
     This class of Calabi-Yau 3-folds occupy an important place in
     mirror symmetry. Now, we want to compute the orbifold
     Dolbeault cohomology of $X$ to compare with Borcea-Voisin's
     calculation of Dolbeault cohomology of $\widetilde{X}$.

         Let's give a brief description of $X$. Our reference is \cite{CB}.
         $\tau$ has 4 fixed points. $(S, \sigma)$ is
         classified by Nikulin. Up to deformation, it is decided
         by three integers $(r, a, \delta)$ with following
         geometric meaning. Let $L^{\sigma}$ be the fixed part of
         $K3$-lattice. Then,
         $$r=rank (L^{\sigma}),
         (L^{\sigma})^*/L^{\sigma}=(\Z/2\Z)^a.\leqno(5.1)$$
         $\delta=0$ if the fixed locus $S_{\sigma}$ of $\sigma$
         represents a class divisible by $2$. Otherwise
         $\delta=1$. There is a detail table for possible value of
         $(r,a,\delta)$ \cite{CB}.

         The cases we are interested in are $(r,a,\delta)\neq (10,
         10,0)$, where $S_{\sigma}\neq \emptyset$. When $(r,a,\delta)\neq
         (10, 8,0)$,
         $$S_{\sigma}=C_g\cup E_1\cdots, \cup E_k \leqno(5.2)$$
         is a disjoint union of a curve $C_g$ of genus
         $$g=\frac{1}{2}(22-r-a)$$
         and $k$ rational curves $E_i$, with
         $$k=\frac{1}{2}(r-a).$$
         For $(r,a,\delta)=(10,8,0)$,
         $$S_{\sigma}=C_1\cup \tilde{C}_1$$
         the disjoint union of two elliptic curves.

         Now, let's compute its orbifold Dolbeault cohomology. We
         assume that $(r, a, \delta)\neq (10, 8, 0)$. The case that
         $(r, a, \delta)=(10, 8, 0)$ can be computed easily as well.
         We leave it as an exercise for the readers.

         An elementary
         computation yields
         $$h^{1,0}(X)=h^{2,0}(X)=0, h^{3,0}(X)=1, h^{1,1}(X)=r+1,
         h^{2,1}(X)=1+(20-r).\leqno(5.3)$$
         Notes that twisted sectors consist of 4 copies of
         $S_{\sigma}$.
         $$h^{0,0}(S_{\sigma})=k+1, h^{1,0}(S_{\sigma})=g.\leqno(5.4)$$
         It is easy to compute that the degree shifting number for
         twisted sectors is 1. Therefore, we obtain
         $$h^{1,0}_{orb}=h^{2,0}_{orb}=0, h^{3,0}_{orb}=1,
         h^{1,1}_{orb}=1+r+4(k+1),
         h^{2,1}_{orb}=1+(20-r)+4g.\leqno(5.5)$$
         Compared with the calculation for $\widetilde{X}$, we get a
         precise agreement.

         Next, we compute the triple product on $H^{1,1}_{orb}$.
         $H^{1,1}_{orb}$ consists of contributions from nontwisted
         sector with dimension $1+r$ and twisted sectors with
         dimension $4(k+1)$. Only nontrivial one is the classes
         from twisted sector. Recall that we need to consider the
         moduli space of 3-point ghost maps with weight $g_1, g_2,
         g_3$ at three marked points satisfying the condition
         $g_1g_2g_3=1$. In our case, the only possibility is
         $g_1=g_2=g_3=\tau\times \sigma$. But $(\tau\times
         \sigma)^3=\tau\times \sigma\neq 1$. Therefore, For any
         class $\alpha$ from twisted sectors, $\alpha^3=0$. On the
         other hand, we know the triple product or exceptional
         divisor of $\widetilde{X}$ is never zero. Hence, $X,
         \widetilde{X}$ have different cohomology ring.
         \vskip 0.1in

         \noindent{\bf Example 5.3-Weighted projective space: } The examples we compute
         so far are global quotient. Weighted projective spaces are the
         easiest examples of non-global quotient orbifolds. Let's consider weighted
         projective space $CP(d_1, d_2)$, where $(d_1, d_2)=1$.
         Thurston's famous tear drop is $CP(1, d)$. $CP(d_1, d_2)$ can be
         defined as the quotient of $S^3$ by $S^1$, where $S^1$
         acts on the unit sphere of $\C^2$ by
         $$e^{i\theta}(z_1, z_2)=(e^{id_1\theta}z_1,
         e^{id_2\theta}z_2).\leqno(5.6)$$
         $CP(d_1, d_2)$ has two singular points $x=[1,0],
         y=[0,1]$. $x,y$ gives rise $d_2-1, d_1-1$ many twisted sectors
         indexed by the elements of isotropy subgroup. The degree
         shifting numbers are $\frac{i}{d_2}, \frac{j}{d_1}$ for
         $1\leq i\leq d_2-1, 1\leq j\leq d_1-1$. Hence, the
         orbifold cohomology are
         $$h^0_{orb}=h^2_{orb}=h^{\frac{2i}{d_2}}_{orb}=
         h^{\frac{2j}{d_1}}_{orb}=1.\leqno(5.7)$$
         Note that orbifold cohomology classes from twisted
         sectors have rational degree. Let $\alpha\in
         H^{\frac{2}{d_1}}_{orb}, \beta\in H^{\frac{2}{d_2}}_{orb}$ be
         the generators corresponding to $1\in H^0(pt, \Q)$. An
         easy computation yields that orbifold cohomology is
         generated by $\{1, \alpha^j, \beta^i\}$ with relation
         $$\alpha^{d_1}=\beta^{d_2},
         \alpha^{d_1+1}=\beta^{d_2+1}=0.\leqno(5.8)$$
         The Poincare pairing is for $1\leq i_1, i_2, i<d_2-1,1\leq j_1, j_2, j<d_1-1$
         $$<\beta^i, \alpha^j>_{orb}=0, <\beta^{i_1},
         \beta^{i_2}>_{orb}=\delta_{i_1, d_2-i_2}, <\alpha^{j_1},
         \alpha^{j_2}>_{orb}=\delta_{j_1, d_1-j_2}.$$
         \vskip 0.1in

    The last two examples are local examples in nature. But they
     exhibit a strong relation with group theory.
    \vskip 0.1in
    \noindent
    {\bf Example 5.4: } The easiest example is probably a point
    with a trivial group action of $G$. In this case, a sector
     $X_{(g)}$ is a point
     with the trivial group action of $C(g)$. Hence,
    orbifold cohomology is
    generated by conjugacy classes of elements of $G$. All the
    degree shifting numbers are zero. Only Poincare pairing and cup
    products are interesting. Poincare paring is obvious. Let's
     consider cup product. First we observe that
    $X_{(g_1,g_2, (g_1g_2)^{-1})}$ is a point with the trivial group
     action of $C(g_1)\cap C(g_2)$. We choose a basis  $\{x_{(g)}\}$
   of the orbifold cohomology group where $x_{(g)}$ is given by the
     constant function $1$ on $X_{(g)}$. Then the inverse of the
     intersection matrix $(<x_{(g_1)},x_{(g_2)}>_{orb})$ has
     $a^{x_{(g)}x_{(g^{-1})}}=|C(g)|$.

Now by Lemma 4.1.4 and Equation $(4.2.12)$, we have
$$
x_{(g_1)}\cup x_{(g_2)}=\sum_{(h_1,h_2), h_1\in (g_1), h_2\in (g_2)}
\frac{|C(h_1h_2)|}{|C(h_1)\cap C(h_2)|} x_{(h_1h_2)},
$$
    where $(h_1,h_2)$ is the conjugacy class of pair $h_1, h_2$.

On the other hand, recall that the center $Z(\C[G])$ of group algebra
    $\C[G]$ is generated by $\sum_{h\in (g)} h$. We can define a
map from the orbifold cohomology group onto $Z(\C[G])$ by
$$
\Psi: x_{(g)}\mapsto \sum_{h\in (g)} h. \leqno (5.9)
$$
The map $\Psi$ is a ring homomorphism, which can be seen as follows:
$$
(\sum_{h\in (g_1)} h)(\sum_{k\in (g_2)}k)=\sum_{h\in (g_1),k\in (g_2)}
hk=\sum_{(h_1,h_2), h_1\in (g_1), h_2\in (g_2)} \frac{A}{B}(\sum_{h\in
(h_1h_2)} h), \leqno (5.10)
$$
where $A=\frac{|G|}{|C(h_1)\cap C(h_2)|}$ is the number of elements in
the orbit of $(h_1,h_2)$ of the action of $G$ given by $g\cdot
(h_1,h_2)= (gh_1g^{-1},gh_2g^{-1})$, and $B=\frac{|G|}{|C(h_1h_2)|}$
is the number of elements in the orbit of $h_1h_2$ of the action of
$G$ given by $g\cdot h=ghg^{-1}$. Therefore, the orbifold cup product
is the same as product of $Z(\C[G])$, and the orbifold cohomology ring
can be identified with the center $Z(\C[G])$ of group algebra
    $\C[G]$ via $(5.9)$.

    \vskip 0.1in
    \noindent
    {\bf Example 5.5: } Suppose that $G\subset SL(n, \C)$ is a
    finite subgroup. Then,
    $\C^n/G$ is an orbifold.  $H^{p, q}(X_{(g)}, \C)=0$ for $p>0$ or $q>0$ and
    $H^{0,0}(X_{(g)}, \C)=\C$.
    Therefore,
    $H^{p, q}_{orb}=0$ for $p\neq q$ and
    $H^{p,p}_{orb}$ is
    a vector space generated by conjugacy class of
     $g$ with $\iota_{(g)}=p$.
    Therefore, we have a natural decomposition
    $$H^*_{orb}(\C^n/G, \C)=Z[\C[G])=\sum_p H_p,\leqno(5.11)$$
    where $H_p$ is generated by conjugacy classes of
    $g$ with $\iota_{(g)}=p$.
    The ring structure is also easy to describe. Let $x_{(g)}$ be
    generator corresponding to zero cohomology class of twisted
    sector $X_{(g)}$. We would
    like to get a formula for $x_{(g_1)}\cup x_{(g_2)}$. As we showed before,
    the multiplication of conjugacy classes can be described in terms of
    center of twisted group algebra $Z(\C[G])$. But we have further
    restrictions in this case. Let's first describe the moduli space
    $X_{(h_1,h_2,(h_1h_2)^{-1})}$ and its corresponding GW-invariants. It is clear
    $$X_{(h_1,h_2,(h_1h_2)^{-1})}=X_{h_1}\cap X_{h_2}/C(h_1, h_2).$$
    To have nonzero invariant, we require that
    $$\iota_{(h_1h_2)}=\iota_{(h_1)}+\iota_{(h_2)}.\leqno(5.12)$$
     Then, we need to compute
   $$\int_{X_{h_1}\cap X_{h_2}/C(h_1,
   h_2)}^{orb} e^*_3(vol_c(X_{h_1h_2}))\wedge e(E),\leqno(5.13)$$
   where $vol_c(X_{h_1h_2})$ is the compact supported $C(h_1h_2)$-invariant top form with volume one
    on $X_{h_1h_2}$. It is also viewed as a form on $X_{h_1}\cap X_{h_2}/C(h_1)\cap
   C(h_2)$.
   However,
   $$X_{h_1}\cap X_{h_2}\subset X_{h_1h_2}$$
   is a submanifold. Therefore, (5.13) is zero unless
   $$X_{h_1}\cap X_{h_2}=X_{h_1h_2}.\leqno(5.14)$$
   In this case, we call $(h_1,h_2)$ transverse. In this case,
   it is clear that obstruction bundle is trivial.  Let
   $$I_{g_1, g_2}=\{(h_1, h_2); h_i \in (g_i),
   \iota_{(h_1)}+\iota_{(h_2)}=\iota_{(h_1h_2)}, (h_1,
   h_2)-transverse \}.\leqno(5.15)$$
   Then, using decomposition lemma 4.1.4,
   $$x_{(g_1)}\cup x_{(g_2)}=\sum_{(h_1,h_2)\in I_{g_1, g_2}}
   d_{(h_1,h_2)} x_{(h_1h_2)}.\leqno(5.16)$$
     A similar computation as previous example yields $d_{(h_1,h_2)}=
     \frac{|C(h_1h_2)|}{|C(h_1)\cap C(h_2)|}$.

\section{Some General Remarks}
         Physics indicated that orbifold quantum cohomology should
         be equivalent to ordinary quantum cohomology of crepant
         resolution.
         It is a rather difficult problem to find the precise relations between orbifold
         quantum cohomology with the quantum cohomology
         of a crepant resolution. At the classical
         level,  there is
         an indication that equivariant $K$-theory is better
         suited for this purpose. For GW-invariant, orbifold GW-invariant
         defined
         in \cite{CR} seems to be equivalent to the relative GW-invariant of
         pairs studied by Li-Ruan \cite{LR}.
         We hope that we will have a better understanding of this relation
         in the near future.

            There are many interesting problems in this orbifold cohomology
            theory. As we mentioned at the beginning, many
            Calabi-Yau 3-folds are constructed as crepant
            resolutions of Calabi-Yau orbifolds. The orbifold string
            theory suggests that there might be a mirror symmetry
            phenomenon for Calabi-Yau orbifolds. Another
            interesting question is the relation between quantum
            cohomology and birational geometry \cite{R},\cite{LR}. In fact, this was
            our original motivation. Namely, we want to
            investigate the change of quantum cohomology under
            birational transformations. Birational transformation
            corresponds to wall crossing phenomenon for symplectic
            quotients. Here, the natural category is symplectic
            orbifolds instead of smooth manifolds. From our work,
            it is clear that we should replace quantum cohomology
            by orbifold quantum cohomology. Then, it is a
            challenge problem to calculate the change of orbifold
            quantum cohomology under birational transformation.
            The first step is to investigate the change of
            orbifold cohomology under birational transformation.
            This should be an interesting problem in its own
            right.

\end{document}